\documentclass[a4paper,reqno,10pt,oneside]{article}
\usepackage{amsmath,geometry}
\usepackage{graphicx}
\usepackage{mathrsfs}
\usepackage{amssymb,amsmath, amsfonts, amsthm}
\usepackage{latexsym}
\usepackage{color} 
\usepackage[dvips]{epsfig}
\usepackage{graphicx}

\allowdisplaybreaks[1]

\numberwithin{equation}{section}

\usepackage{pdfsync}

\title{Principal eigenvalues and eigenfunctions for fully nonlinear equations in punctured balls}

\author{Isabeau Birindelli \\
Dipartimento di Matematica, Sapienza Universit\`a\  di Roma
\and
  Fran\c{c}oise Demengel\\
  D\'epartement de Math\'ematiques,
CY Paris University
  \and
   Fabiana  Leoni\\ 
   Dipartimento di Matematica, Sapienza Universit\`a\  di Roma}

\date{}

\usepackage{amsthm}
\usepackage{amsmath}

\catcode`@=11 \@addtoreset{equation}{section} \catcode`@=12

\newtheorem{theo}{Theorem}[section]
\newtheorem{prop}[theo]{Proposition}
\newtheorem{rema}[theo]{Remark}

\newtheorem{cor}[theo]{Corollary}
\newtheorem{lemme}[theo]{Lemma}

\def\R{\mathbb  R}

\newcommand{\N}{{\mathbb N}}

\begin{document}
\maketitle
\begin{abstract} This paper is devoted to the proof of the existence of  the principal eigenvalue and related eigenfunctions
for fully nonlinear uniformly elliptic equations posed in a punctured ball, in presence of  a singular potential. More precisely,   we   analyze  existence, uniqueness  and regularity    of solutions $( \bar\lambda_\gamma, u_\gamma)$ of the equation 
$$F( D^2 u_\gamma)+ \bar \lambda_\gamma \frac{u_\gamma}{r^\gamma} = 0\ {\rm in} \ B(0,1)\setminus \{0\}, \ u_\gamma = 0 \ {\rm on} \ \partial B(0,1)$$
 where $u_\gamma>0$ in $B(0,1)$, and $\gamma >0$. We  prove existence of radial solutions which are continuous on $\overline{ B(0,1)}$ in the case $\gamma <2$,  existence of unbounded solutions in the case $\gamma = 2$ and a non existence result for $\gamma >2$. 
 We also give the explicit value of $\bar \lambda_2$ in the case of  Pucci's operators,  which generalizes the  Hardy--Sobolev constant for the Laplacian. 
\end{abstract}
\section{Introduction}

In this paper we will study radial eigenvalues  and related positive radial eigenfunctions for the  Dirichlet problem 
  $$F( D^2 u) +\mu  r^{-\gamma} u = 0 \,  \mbox{ in }\ \overline{B(0,1)} \setminus \{0\}$$
 when  $\gamma>0$, and $F$ is a second order fully nonlinear uniformly elliptic operator. By radial eigenvalue and radial eigenfunction we mean respectfully a real value $\lambda_\gamma$ and a radial nontrivial function  $u_\gamma$   satisfying  the equation 
 \begin{equation}\label{P}
 \left\{\begin{array}{lc}
 F( D^2 u_\gamma) +\lambda_\gamma  r^{-\gamma} u_\gamma = 0  & \mbox{ in }\ \overline{B(0,1)} \setminus \{0\}\\
 u_\gamma=0 & \mbox{ on }\ \partial B(0,1) .
 \end{array}
 \right.
 \end{equation}
 In principle, eigenfunctions are required to satisfy the above eigenvalue problem in the viscosity sense, but, due to the radial symmetry, this is equivalent to consider  classical solutions.

  We will focus on constant sign eigenfunctions, in particular positive eigenfunctions, thus  referring  to the so called  principal eigenvalues.  If necessary, in order to emphasize the dependence of the eigenvalue on the operator $F$, the potential $f(r)$ appearing in the zero order term and the domain $\Omega$ in which the equation is considered, we   will use the notation $\lambda = \lambda(F, f(r), \Omega)$.
  
Interestingly,  we will see that for problem \eqref{P}, as in the case when $F$ is the  Laplace operator,  $\gamma = 2$  is a critical value, in the sense that for $\gamma <2$  there exists smooth eigenfunctions, for $\gamma >2$ there are no eigenfunctions and,   for $\gamma=2$, the eigenfunctions are unbounded.
 \smallskip
 
 Let us recall some  known results  when $F$ is the  Laplacian.
 In the case $\gamma = 2$, the equation is naturally linked  to Hardy's inequality.
Indeed, if $N>2$ and  $u\in H_0^1 ( B(0,1))$ (respectively, $u\in H^1( \R^N)$), then $\frac{u(x)}{|x|}$ belongs to $L^2(B(0,1))$ (respectively $\frac{u(x)}{|x|} \in L^2( \R^N)$), and there exists a positive constant $c$ such that 
  $$ \int\left(  \frac{|u(x)|}{|x|}\right)^2 \leq c \int | \nabla u |^2.$$
 Furthermore, the best constant  $c= \frac{4}{(N-2)^2}$  is not achieved, in the sense that 
  \begin{equation}\label{varlap} \inf _{ u\in H_0^1( B(0,1)), \int_{B(0,1)} \left(\frac{|u(x)|}{|x|}|\right)^2 =1}  \int_{B(0,1)}  | \nabla u |^2= \frac{(N-2)^2}{4} \end{equation}
but there is no $u\in H_0^1$ which realizes the infimum.   By obvious arguments, $B(0,1)$ can be replaced by any bounded regular open set of $\R^N$ containing $0$,  and the optimal constant does not depend on the size of $\Omega$. Note that the right hand side of (\ref{varlap}) coincides with the variational characterization of  the first (or principal) eigenvalue for the equation 
   $$-\Delta u = \lambda \frac{u}{|x|^2}.$$
   For further knowledge on the Hardy--Sobolev inequality and for the case of the $p$-Laplacian, we refer to  \cite{GAPA, PDe, St}.

   On the other hand, if the exponent $\gamma $ of the potential  is strictly less than 2, since  $H_0^1(B(0,1))$ is compactly embedded into the weighted space $L^2 ( B(0,1), \frac{1}{r^\gamma})$, 
  then  existence of minima in  $H_0^1 ( B(0,1))$ can be obtained by  standard arguments of  the direct method in calculus of variations. In that case, denoting 
 \begin{equation}\label{variational}
 \bar \lambda_\gamma =  \inf _{ u\in H_0^1( B(0,1)), \int_{B(0,1)}\frac{|u(x)|^2}{|x|^\gamma} =1}  \int _{B(0,1)}| \nabla u |^2
 \end{equation}
    one sees that $\bar \lambda_\gamma$ is also the first eigenvalue for the equation
    $$\Delta u + \bar \lambda _\gamma \frac{u}{r^\gamma}  = 0\, ,$$ meaning  that $\bar \lambda_\gamma$ is such that there exists $u>0$ in $H_0^1( B(0,1))$ satisfying the equation. 
    
     Note that, by its definition, $\bar \lambda_\gamma$ depends on the domain, since
     $$\bar \lambda_\gamma ( B(0, t) ) = \frac{1}{t^{ 2-\gamma}}\bar \lambda _\gamma ( B(0,1)).$$
     
      If $\gamma >2$ there is no embedding from  $H_0^1(B(0,1))$ into $L^2 ( B(0,1), \frac{1}{r^\gamma})$. Indeed, as an example, the function $$u(r) = r^{-\frac{N-2}{2}+\epsilon}(-\log r)$$ with $0<\epsilon < \frac{\gamma-2}{2}$, belongs to $H_0^1 ( B(0,1))$ and satisfies 
      $$ \int _{ B(0,1) }\frac{u(|x|)^2}{|x|^\gamma} = +\infty$$
 For results in the variational linear case we refer to the works of  many authors, but in particular we wish to mention the works of Cirstea and collaborators \cite{Ci, CDu1, CDu2, CirsChau} .
      
Let us  now focus on  the case of concern of this paper i.e. when $F$ is a fully nonlinear uniformly elliptic operator, that is $F$ is a continuous function defined on the set $\mathcal{S}_N$ of symmetric $N\times N$ matrices, and it satisfies, for positive constants $\Lambda\geq \lambda>0$,
       \begin{equation}\label{FNL} \lambda\,  \hbox{tr}( M^\prime) \leq F( M+ M^\prime) -F( M) \leq \Lambda\,  \hbox{tr}(M^\prime)\, ,
       \end{equation}
for all $M, M'\in \mathcal{S}_N$, with $M'$ positive semidefinite.

We suppose also that $F$ is rotationally invariant, that is
\begin{equation}\label{rotin}
F(O^t MO)= F(M)
\end{equation}
for every orthogonal matrix $O$ and for all $M\in \mathcal{S}_N$, and that $F$ is  positively homogeneous of degree $1$, i.e. 
       \begin{equation}\label{posh} 
       F( tM) = tF(M)\end{equation}
for any $M\in \mathcal{S}_N$ and for all $t>0$. In this case we will see that, as in the regular case i.e. $\gamma=0$, the first eigenvalue for problem \eqref{P}
can be defined  on the model  of \cite{BNV}, i.e. by the optimization formula 
        \begin{equation}\label{deflambdagamma} 
\begin{array}{rl}
        \bar \lambda_\gamma \!\!\! & = \bar \lambda_\gamma (F, r^{-\gamma}, B(0,1)\setminus \{0\})\\[2ex]
        & = \sup \{\mu\, : \  \exists\,   
        u \in C(B(0,1)\setminus \{0\})\, ,\ u>0 \hbox{ in } B(0,1)\setminus \{0\},   \ F (D^2 u) + \mu \frac{u}{r^\gamma} \leq 0\}\, ,
        \end{array}
        \end{equation} 
where the differential inequality  is understood in the viscosity sense.      
       
       The first easy observation is that, by considering constant super-solutions,  one always has $\bar \lambda_\gamma\geq 0$. 
One of the goals  of the present paper is to show,  in particular, that $\bar \lambda_\gamma>0$ for  $\gamma\leq 2$. 

In case of smooth coefficients and regular domains, the principal eigenvalues and related eigenfunctions    
for fully  nonlinear operators $F$ have been largely investigated. We refer to e.g. \cite{IY, BD1, BD2, BLP, QS}.

The  Pucci's extremal operators will play a crucial role and we will treat them in depth. We begin by recalling their definition: by decomposing each matrix $M\in \mathcal{S}_N$ as $M = M^+-M^-$, where $M^+$ and $M^-$ are positive semidefinite matrices satisfying  $M^+M^-=O$, then Pucci's sup operator can be defined as
$$\mathcal{M}^+ _{\lambda, \Lambda} (M) = \Lambda \hbox {tr}( M^+)- \lambda \hbox {tr}( M^-)\, ,$$
 as well as Pucci's inf operator is given by
           $$\mathcal{M}^-_{\lambda, \Lambda} (M)= \lambda \hbox {tr}( M^+)- \Lambda \hbox {tr}( M^-)= -\mathcal{M}^+_{\lambda, \Lambda}  (-M).$$
 As it is well known, see \cite{CC}, under assumptions  \eqref{FNL} and \eqref{posh}, each operator $F$ satisfies
 \begin{equation}\label{elliptic}
 \mathcal{M}^-_{\lambda ,\Lambda}(M)\leq F(M)\leq \mathcal{M}^+_{\lambda ,\Lambda}(M)\, ,\qquad \forall\, M\in \mathcal{S}_N\, ,
 \end{equation}
 showing as Pucci's operators act as explicit extremal   operators    in the whole class of uniformly elliptic operators having the same ellipticity constants. In the sequel, we will omit in the notation the dependence on  the ellipticity constants, which are fixed once for all.    
 
 We further recall that for a   $C^2$  radial function $u(x) = u(|x|)$, one has
                      $$ D^2 u(x)  = u^{\prime \prime} (r)\frac{x\otimes x }{r^2} + \frac{u^\prime (r)}{r} \left( I-\frac{x\otimes x}{r^2}\right)\, ,$$
and, as a consequence,  
            $${\cal M}^+ ( D^2 u) =\Lambda (N-1)\left(\frac{u^\prime (r)}{r}\right)^+- \lambda (N-1)\left(\frac{u^\prime (r)}{r}\right)^-+ \Lambda (u^{\prime \prime}(r))^+-\lambda  ( u^{\prime \prime}(r) )^-\, ,$$
 $${\cal M}^- (D^2 u) =\lambda (N-1)\left(\frac{u^\prime (r)}{r}\right)^+- \Lambda (N-1)\left(\frac{u^\prime (r)}{r}\right)^-+ \lambda (u^{\prime \prime}(r))^+-\Lambda  ( u^{\prime \prime}(r) )^-\, .$$
 Thus,  the ODEs satisfied by radial solutions  of Pucci's extremal equations have coefficients depending on the  dimension like parameters, associated with $\mathcal{M}^+$ and $\mathcal{M}^-$ respectively, defined as     
$$\tilde N_+ = \frac{\lambda}{\Lambda} ( N-1) +1\, ,\quad \tilde N_- = \frac{\Lambda}{\lambda}( N-1)+1\, .$$
 Note that one has always 
$$
\tilde{N}_-\geq N \geq \tilde{N}_+\, ,
$$
with equalities holding true if and only if $\Lambda=\lambda$. We will  assume always that $\tilde{N}_+>2$.

We can now state the main results  of the paper. We will always assume that $F$ satisfies assumptions \eqref{FNL}, \eqref{rotin} and \eqref{posh}. Let us start with the case $\gamma<2$. 
\begin{theo} \label{exigamma}
Suppose that $\gamma <2$.    Then:
 \begin{itemize}   
\item[(i)] $\bar \lambda_\gamma $ defined in (\ref{deflambdagamma} ) is positive and there exists a  function $u$, continuous in  $\overline{B(0,1)}$, radial, strictly positive in $B(0,1)$, such that 
     $$\left\{ \begin{array}{cl}
     F ( D^2 u) + \bar \lambda_\gamma \frac{u}{r^\gamma} = 0 & \hbox{ in } \ B(0,1)\setminus \{0\}\\[1ex]
       u = 0 & \hbox{ on } \ \partial B(0,1)
       \end{array}\right..$$
       Furthermore $u$ is $C^2(B(0,1)\setminus\{0\})$ and it can be extended on  $B(0,1)$ as a Lipschitz continuous function if $\gamma \leq 1$, as a function of class $C^1(B(0,1))$ when $\gamma <1$,  and as an H\"older continuous function with exponent $2-\gamma$ if $\gamma >1$. 
       
\item[(ii)] $\bar \lambda_\gamma$ is stable under various regular approximations : 
 $$\bar \lambda_\gamma = \lim_{ \epsilon \rightarrow 0} \bar \lambda (F,   \frac{1}{(r^2+\epsilon^2)^{\frac{\gamma}{2}}}, B(0,1))\, ,\,$$
      $$\bar \lambda_\gamma = \lim_{ \delta \rightarrow 0 } \bar \lambda (F, \frac{1}{r^\gamma},  B(0,1) \setminus \overline{ B(0, \delta)})\, .$$
  \end{itemize}       
\end{theo}
      
Statement (i) of the  above theorem shows in particular that $\bar \lambda_\gamma$ is actually achieved on smooth radial eigenfunctions. Thus,  
      if we define 
\begin{equation}\label{smoothlambda}  \bar \lambda_\gamma^{\prime} := \sup \{\mu\, : \  \exists\,   
        u \in C^2(B(0,1)\setminus \{0\})\, ,\ u>0 \hbox{ in } B(0,1)\setminus \{0\}, \ u \hbox{ radial},  \ \ F (D^2 u) + \mu \frac{u}{r^\gamma} \leq 0\}
 \, ,
 \end{equation}  
  it then follows that 
  $$\bar \lambda_\gamma = \bar \lambda_\gamma^\prime.$$
 Actually, we will work initially with the smooth eigenvalue $\bar \lambda_\gamma'$, and we will finally prove that it coincides with $\bar \lambda_\gamma$.  Note that, due to the lack of regularity of the coefficient function $\frac{1}{r^\gamma}$ we cannot employ directly the results of \cite{DaLS}  which ensure that solutions of 
   $$ F( D^2 u) + f(r) u=0$$
    in a radial domain  are radial when $f$ is non increasing.   Nonetheless,  we will prove   that the two eigenvalues coincide for any $\gamma\leq 2$. However, due to the singularity at zero,  we cannot prove that  any eigenfunction is radial.  
    
    Theorem \ref{exigamma} will be proved after several steps and intermediate results. In particular, we will prove a comparison theorem for smooth, bounded, radial sub- and super-solutions in the punctured ball, without assuming any order condition at the origin. Furthermore, we will show that for $\mu<\bar \lambda_\gamma$ the problem  
$$
 \left\{
 \begin{array}{cl} 
  F( D^2 u)+ \mu u  r^{-\gamma}  =  f (r) r^{-\gamma} & \hbox{ in } B(0,1)\setminus \{0\}\\
   u=0 & \hbox{ on } \partial B(0,1)
 \end{array}\right.  $$
admits a unique radial solution $u\in C^2(B(0,1)\setminus \{0\}) \cap C(\overline{B(0,1)})$ for any  radial and continuous datum $f\in \overline{B(0,1)}$ satisfying $f\leq 0$. 

Next, for the case $\gamma=2$ we have the following result, which gives explicit expressions for the eigenvalues and the eigenfunctions in case of Pucci's operators.
 
       \begin{theo}\label{gamma=2}
 Assume that  $\gamma = 2$. Then: 
\begin{itemize}
\item[(i)]   For the operator $\mathcal{M}^+$ one has
$$ \bar \lambda_2(\mathcal{M}^+) =\Lambda \frac{(\tilde N_+-2)^2}{4}$$
and the function
$
u(x)= r^{-\frac{\tilde{N}_+-2}{2}} (-\ln r)
$
is an explicit solution of
$$\left\{ \begin{array}{cl}
 \mathcal{M}^+( D^2 u) + \bar \lambda_2 \frac{u}{r^2} = 0 & \hbox{ in } \ B(0,1)\setminus \{0\}\\[1ex]
       u = 0 & \hbox{ on } \ \partial B(0,1)
       \end{array}\right.$$
   Analogously, for the operator $\mathcal{M}^-$ one has
$$ \bar \lambda_2(\mathcal{M}^-) =\lambda \frac{(\tilde N_--2)^2}{4}$$
and the function
$
u(x)= r^{-\frac{\tilde{N}_--2}{2}} (-\ln r)
$
is an explicit solution of
$$\left\{ \begin{array}{cl}
 \mathcal{M}^-( D^2 u) + \bar \lambda_2 \frac{u}{r^2} = 0 & \hbox{ in } \ B(0,1)\setminus \{0\}\\[1ex]
       u = 0 & \hbox{ on } \ \partial B(0,1)
       \end{array}\right.$$
\item[(ii)]  The eigenvalues $\bar \lambda_2(\mathcal{M}^\pm)$ are stable under various  regularization 
     $$\bar \lambda_2(\mathcal{M}^\pm) = \lim_{ \gamma \rightarrow 2} \bar \lambda_\gamma(\mathcal{M}^\pm)$$ 
      $$\bar \lambda_2(\mathcal{M}^\pm) = \lim_{ \delta \rightarrow 0 } \bar \lambda (\mathcal{M}^\pm,  B(0,1) \setminus \overline{ B(0, \delta)})$$
      $$\bar \lambda_2(\mathcal{M}^\pm) = \lim_{ \epsilon \rightarrow 0} \bar \lambda  (\mathcal{M}^\pm, \frac{1}{( r^2+\epsilon^2)})$$
\item[(iii)] For any operator $F$ satisfying \eqref{FNL} and \eqref{posh} one has
$$
\Lambda \frac{(\tilde N_+-2)^2}{4}\leq \bar \lambda_2 (F)\leq \lambda \frac{(\tilde N_--2)^2}{4}\, .
$$     
   \end{itemize}
    \end{theo}
We observe that we cannot prove the existence of eigenfunctions for a general operator $F$, but we can merely provide the estimate on the eigenvalue given by statement (iii) above.    
Theorem \ref{gamma=2} will be obtained by using a variational approach adapted to the fully nonlinear radial framework. Indeed, we will define variational eigenvalues associated with the operators $\mathcal{M}^\pm$ in an analogous way as in \eqref{variational}, taking advantage of the radial symmetry of solutions. Then, the full statements of Theorem \ref{gamma=2} will follow as consequences of the properties established for $\bar \lambda_\gamma$ in the case  $\gamma<2$ and the stability of the variational formulation as $\gamma\to 2$.
    
Finally, for the case $\gamma>2$, the singularity of the coefficient is too strong and it prevents the existence of positive smooth super-solutions, as stated by the following non existence result.   
    \begin{theo}\label{gamma>2}
  If $\gamma >2$, then  the eigenvalue $\bar \lambda^\prime_\gamma$ defined by \eqref{smoothlambda} satisfies    $\bar \lambda_\gamma^\prime = 0$.       \end{theo}
    
 Let us observe that,
  symmetrically, one could  define the eigenvalue associated with negative eigenfunctions, by setting
  $$ \bar \lambda_\gamma^- =  \sup \{\mu\, : \  \exists\,   
        u \in C(B(0,1)\setminus \{0\})\, ,\ u<0 \hbox{ in } B(0,1), \  \ \ F (D^2 u) + \mu  \frac{u}{r^\gamma} \geq 0\}=0
   . $$ 
In this case,  the results above can be extended to $\bar \lambda_\gamma^- $ with obvious modifications.

\bigskip

Let us conclude this introduction by observing that, in case of semilinear or quasilinear equations, many existence, non existence and classification results have been obtained in presence  of zero oder terms having  Hardy's potential perturbed with additional sub- or superlinear terms. In particular, we refer to \cite{BoO, BOP,  Ci, W Du} for results related to Laplace operator, and to \cite{GAPA} for the $p$-Laplace operator. 
      
The case where, in all  directions above, $\Delta$ or $\Delta_p$ is replaced by a non variational  fully nonlinear operator   will be  the object of  future works. 
  

  \section{ The case  $ \gamma <2$: proof of Theorem \ref{exigamma} }
Theorem \ref{exigamma} will be proved as a consequence of  several classical steps: a comparison principle,  existence and regularity results  and a maximum principle "below" the first eigenvalue. 


\subsection{Maximum principles, existence and regularity results}

 The first result of the present section is a crucial technical lemma. 
 \begin{lemme}\label{fabiana}
Let $ f\in C\left( B(0,1)\setminus \{0\}\right)$ be a radial, bounded and  positive function and assume that  $u\in C^2\left( B(0,1)\setminus \{0\}\right)$ is  a radial,  bounded function satisfying
\begin{equation}\label{subsol}{\cal M}^+ ( D^2 u) \geq  f   r^{-\gamma}\qquad \hbox{ in } B(0,1)\setminus \{0\}\, .\end{equation}
 Then 
\begin{itemize}
 \item[(i)]  $u^\prime \geq 0$ in a right neighborhood of $0$;
  \item[(ii)] $\displaystyle \lim_{r\rightarrow 0} u^\prime(r) r^{ {\tilde N}_--1} = 0$ and in a right neighborhood of $0$ one has
  $$ u'(r)\geq    \frac{\inf f}{\Lambda ({\tilde N}_--\gamma)} r^{1-\gamma}\, ;$$
 \item[(iii)] if,  furthermore, ${ \cal M}^+ ( D^2 u) = f r^{-\gamma}$, then, in a   right neighborhood of $0$, one has also
                      $$u'(r) \leq \frac{\sup f}{\lambda (N-\gamma)} r^{1-\gamma}\, .$$
In particular, there exists a constant $c>0$ such that, for $r$ sufficiently small,
                       $$|u^\prime (r) | \leq c r^{1-\gamma}$$
                        and then $u$ is locally Lipschitz continuous in $B(0,1)$ if $\gamma \leq 1$,  it belongs to  $C^1(B(0,1))$ if $\gamma <1$, and it is locally H\"older continuous in $B(0,1)$ with exponent $2-\gamma$ if $\gamma >1$. 
 \end{itemize}
                        \end{lemme} 
                                               \begin{proof}
 Let us prove, by contradiction, that   $u^\prime $ does not change sign in a right neighborhood of 0. If not, there exists a decreasing sequence $\{r_n\}$ converging to 0, such that $u'(r_n)=0$ for all $n$, and  $u^\prime \leq 0$ in $]r_{2n+1}, r_{2n}[$,    $u^\prime \geq 0$ in $]r_{2n+2}, r_{2n+1}[$. Since $u^\prime (r_{2n}) = u^\prime (r_{2n+1}) = 0$, there exists some $s_{2n}\in ]r_{2n+1}, r_{2n}[$ such that $u^{\prime \prime} ( s_{2n}) = 0$. This yields the contradiction
 $$0\geq \lambda (N-1) \frac{u'(s_{2n})}{s_{2n}}=\mathcal{M}^+(D^2u(s_{2n}))\geq f(s_{2n}) s_{2n}^{-\gamma}>0\, .
 $$
Next, arguing again by contradiction, if  $u'(r)\leq 0$ for $r$ sufficiently small, then, by \eqref{subsol}, one has
                  $u^{\prime \prime}(r) >0$ and, in a right neighborhood of $0$,
$$ \mathcal{M}^+(D^2u)=\Lambda u''(r) +\lambda (N-1) \frac{u'(r)}{r} \geq f(r) r^{-\gamma}\, .$$
Hence,
$$
(u^\prime(r) r^{ {\tilde N}_+-1})^\prime  \geq \frac{f  r^{ {\tilde N}_+-1-\gamma}}{\Lambda}>0$$
and,  in particular, $u^\prime  r^{{\tilde N}_+-1}$ is increasing in a right neighborhood of $0$. Thus,  $\lim_{r\to 0}
u^\prime(r)  r^{ {\tilde N}_+-1}$   exists and it is lesser than or equal to 0. If it was lesser than zero, then we would have, for   some constant $l>0$, 
$$ u^\prime(r) \leq - l r^{1- {\tilde N}_+} $$
in a right neighborhood of  $0$, yielding a contradiction to  the boundedness of  $u$   This shows that $\lim_{r\to 0}
u^\prime(r)  r^{ {\tilde N}_+-1}=0$ and, by monotonicity, $u'(r)>0$ in a right neighborhood of $0$. The reached contradiction proves statement (i).
                     
In order to prove (ii), let us observe that, for $r$ sufficiently small,  by \eqref{subsol} we have either                        
 $$  u''(r) +(N-1) \frac{u'(r)}{r} \geq \frac{f(r)r^{-\gamma}}{\Lambda}$$
 if $u''(r)\geq 0$, or
$$ u''(r) +({\tilde N}_--1) \frac{u'(r)}{r} \geq \frac{f(r)r^{-\gamma}}{\lambda}$$
if $u^{\prime \prime}(r) \leq 0$. Since $u'(r)\geq 0$ and ${\tilde N}_-\geq N$, in both cases one has
$$ u''(r) +({\tilde N}_--1) \frac{u'(r)}{r} \geq \frac{f(r)r^{-\gamma}}{\Lambda}\, ,$$
that is                    
 $$ (u^\prime r^{{\tilde N}_--1} )^\prime \geq \frac{f(r) r^{{\tilde N}_--1-\gamma}}{\Lambda}\geq \frac{\inf f}{\Lambda} r^{{\tilde N}_--1-\gamma}\, .$$
 Arguing as above, we deduce that
                        $ u^\prime r^{{\tilde N}_--1}$ is increasing in a right neighborhood of $0$, hence it has a nonnegative limit as $r\to 0$, and  such a limit must be  $0$, since $u$ is bounded.  Moreover, by integrating the above inequality, we obtain
 $$
 u'(r)\geq    \frac{\inf f}{\Lambda ({\tilde N}_--\gamma)} r^{1-\gamma}\, .
 $$                 
Let us finally  prove (iii). Assuming that $\mathcal{M}^+(D^2u) =f(r)r^{-\gamma}$ and using statement (i), it follows that, for every $r>0$ sufficiently small, one has either
$$  u''(r) +(N-1) \frac{u'(r)}{r} = \frac{f(r)r^{-\gamma}}{\Lambda}$$
or
$$  u''(r) +({\tilde N}_--1) \frac{u'(r)}{r} = \frac{f(r)r^{-\gamma}}{\lambda}\, .$$
In both cases, we deduce
$$u''(r) +(N-1) \frac{u'(r)}{r} \leq \frac{f(r)r^{-\gamma}}{\lambda}\, ,$$
which yields
$$ (u^\prime(r) r^{ N-1})^\prime \leq \frac{f(r)}{\lambda}  r^{N-1-\gamma}\leq \frac{\sup f}{\lambda} r^{N-1-\gamma}\, . $$
Hence, $u'(r) r^{N-1}-\frac{\sup f}{\lambda (N-\gamma)} r^{N-\gamma}$ is non increasing in a right neighborhood of $0$ and it has a limit as $r\to 0$. This implies that  $u'(r) r^{N-1}$ has a limit as $r\to 0$ as well, and such a limit is zero by the boundedness of $u$.  By integrating the last inequality, we  finally deduce
$$
u'(r) \leq \frac{\sup f}{\lambda (N-\gamma)} r^{1-\gamma}\, .
$$
The regularity of $u$ at zero is then a consequence of the estimate $|u'(r)|\leq c r^{1-\gamma}$. Elsewhere, it follows from the assumption $u\in C^2(B(0,1)\setminus \{0\})$.

 \end{proof}
 
  \begin{rema}\label{lemmaformm}
 {\rm  By using the change of variable $v = -u$, one gets that if $ f\in C\left( B(0,1)\setminus \{0\}\right)$ is  a  radial, bounded and  positive function and   $u\in C^2\left( B(0,1)\setminus \{0\}\right)$ is  a  bounded radial function satisfying
$$ { \cal M}^-(D^2 u) \leq -f r^{-\gamma}\quad \hbox{ in } B(0,1)\setminus \{0\}\, ,$$
  then, for $r$ sufficiently small, $u^\prime (r) \leq 0$, 
                         $\lim_{r\to 0} u^\prime (r)r^{{\tilde N_--1}}=0$ and
$$ 
u'(r)\leq    -\frac{\inf f}{\Lambda ({\tilde N}_--\gamma)} r^{1-\gamma}\, .
$$
Moreover, if  ${ \cal M}^-(D^2 u) = -f r^{-\gamma}\quad \hbox{ in } B(0,1)\setminus \{0\}$, then $|u^\prime (r) | \leq c r^{1-\gamma}$ for a positive constant $c$. Hence  $u$ is locally Lipschitz continuous in $B(0,1)$ for $\gamma\leq 1$, it belongs to $C^1(B(0,1))$ if $\gamma <1$,  and it is locally H\"older continuous in $B(0,1)$ with exponent $2-\gamma$ for $\gamma >1$. 
 
Obviously, since 
                          ${\cal M}^-\leq { \cal M}^+$, one gets an analogous  conclusion when 
                          $${ \cal M}^+ ( D^2 u) \leq -f r^{-\gamma}.$$}
 \end{rema} 
 
 We can now prove a comparison principle for general radial fully nonlinear singular equations of the form
 $$
 F(D^2u)- \beta u r^{-\gamma} = f(r) r^{-\gamma} \quad \hbox{ in } B(0,1)\setminus \{0\}
 $$
 when no boundary condition at the origin is assumed.

\begin{theo}\label{compa}
Let $f, g\in C\left(B(0,1)\right)$ be  radial functions and assume that $u, v\in C\left( \overline{B(0,1)}\right) \cap C^2(B(0,1)\setminus \{0\})$ are radial functions satisfying in $B(0,1)\setminus \{0\}$
 $$
 \begin{array}{c}                 
F( D^2 u) -\beta u(r) r^{-\gamma} \geq f(r)r^{-\gamma} \\[2ex]
F( D^2 v) -\beta v(r) r^{-\gamma}  \leq g(r) r^{-\gamma}
\end{array}
$$
with $\beta\geq 0$ and   $f\geq g$ in $B(0,1)$.  Then,  $u\leq v$ on $\partial B(0,1)$ implies $u\leq v$ in 
$\overline{B(0,1)}$.
 \end{theo}
              
\begin{proof} 

Let us first consider the case in which either $\beta>0$ or  $f>g$ in $B(0,1)$. We           
suppose by contradiction that 
$$\max_{\overline{B(0,1)}} (u-v) >0\, .
$$
If the maximum is achieved at 0, then $(u-v)(0) >0$ and, by the assumptions on $f$, $g$ and $\beta$, there exist $\delta>0$ and  a neighborhood on the right of $0$ on which 
               $$ (\beta (u(r)-v(r) )+f(r)-g(r))r^{-\gamma} \geq \delta r^{-\gamma} >0$$
By the uniform ellipticity of $F$, we then obtain for $r$ sufficiently small
$${ \cal M}^+ (D^2(u-v) )\geq \delta r^{-\gamma}$$
Using Lemma \ref{fabiana} for $u-v$, one gets that  for some positive constant $c$, $(u-v)^\prime \geq  c \delta r^{1-\gamma}$, which contradicts the fact that   $u-v$ attains its maximum  at $0$. 
Hence, there exists $0<\bar r<1$ such that  $u(\bar r)-v(\bar r)=\max (u-v)$. Then,   $(D^2 u-D^2 v)(\bar r) \leq 0$ and, by ellipticity, we get   the contradiction
                    $$  (f( \bar r) + \beta u(\bar r))r^{-\gamma}  \leq F( D^2 u(\bar r) )\leq F( D^2 v (\bar r))\leq   (g(\bar r) + \beta v(\bar r)) r^{-\gamma}\, .$$
For the case $\beta=0$ and $f\geq g$, let us introduce the radial function
$$
w(r) =  1-r^\tau$$
with $0< \tau \leq 2-\gamma$. A direct computation shows that
$$
\mathcal{M}^+(D^2w)\leq \tau \Lambda ( |\tau-1| - (\tilde{N}_+-1) ) r^{\tau-2}\, .
$$
We observe that $|\tau -1|<1<\tilde{N}_+-1$, so that
$$
\mathcal{M}^+(D^2w)\leq -Cr^{-\gamma}\qquad \hbox{ in } B(0,1)\setminus \{0\}
$$
with  $C=\tau \Lambda (\tilde{N}_+-1-|\tau -1|)>0$. Thus, for any $\epsilon>0$, we have
$$
F(D^2(u-\epsilon w))\geq F(D^2u)- \epsilon \, \mathcal{M}^+(D^2w)\geq (f+\epsilon\, C)r^{-\gamma}\, .
$$
Since $f+\epsilon\, C>g$ and $u-\epsilon w=u\leq v$ on $\partial B(0,1)$, the previous argument proves that
$$
u-\epsilon w\leq  v \qquad \hbox{ in } \overline{B(0,1)}
$$
and the conclusion follows by letting $\epsilon\to 0$.

 \end{proof}   
 
 \begin{rema}\label{posi}
{\rm      The auxiliary function  introduced in the proof of Theorem \ref{compa} shows that there exist a  radial function $w\in C^2(B(0,1)\setminus \{0\})$, strictly positive in $B(0,1)\setminus \{0\}$,  
such  that 
                      $$ F( D^2 w) \leq
                       - cw  r^{-\gamma}\qquad \hbox{ in } B(0,1)\setminus \{0\}$$ 
for a constant $c>0$. This proves that                      
 $\bar \lambda_\gamma (F)\geq \bar \lambda_\gamma' (F)\geq c>0$. }
 \end{rema}              

Next, we have the following existence, uniqueness and regularity  result.
                              
\begin{theo}\label{exi1}
              Let $f \in C (B(0,1))$ be a radial, bounded function. For $\beta \geq 0$ and $b\in \R$ 
   there exists a unique bounded radial function $u\in C(\overline {B(0,1)}\setminus \{0\})\cap C^2(B(0,1)\setminus \{0\})$  satisfying
 \begin{equation}\label{dp}
 \left\{\begin{array}{cc}
               F( D^2 u) -\beta{ u r^{-\gamma}}  =  r^{-\gamma}   f(r) & \hbox{ in } \ B(0,1)\setminus \{0\}\\
                u = b & \hbox{ on } \ \partial B(0,1) 
                \end{array}\right.
 \end{equation}               
Moreover, $u$ can be extended up to $\overline{B(0,1)}$, and one has: $u\in C^1(\overline{B(0,1)})$ if $\gamma<1$, $u$ is Lipschitz continuous in $\overline{B(0,1)}$ if $\gamma \leq 1$, $u$ is H\"older continuous in  $\overline{B(0,1)}$ with exponent $2-\gamma$ if $\gamma > 1$.            
              \end{theo}
 \begin{proof} 
 For every $n\in \N$ let us introduce the regularized Dirichlet boundary value problem
 \begin{equation}\label{dpapp}
 \left\{\begin{array}{cl}
 F( D^2 u_n) -\beta{ u_n (r^2+1/n)^{-\gamma/2}}  =  (r^2+1/n)^{-\gamma/2}   f(r) & \hbox{ in } \ B(0,1)\\
                u_n = b & \hbox{ on } \ \partial B(0,1) 
                \end{array}\right.
 \end{equation}  
 which, by standard viscosity solutions theory, see \cite{usr}, has a unique solution $u_n\in C(\overline{B(0,1)})$. By the symmetry results of \cite{DaLS}, it follows that $u_n$ is radial, hence, as a solution of the associated ODE,  $u_n$ belongs to $C^2(B(0,1))$.

For $0<\tau\leq 2-\gamma$, let us consider the radial function
 $$
 w(r)=L(1-r^\tau) +b\, ,
 $$
 where $L>0$ is a constant to be suitably chosen.
 The same computation used in the proof of Theorem \ref{compa} yields
 $$
 \mathcal{M}^+(D^2w)\leq -L\, C  r^{-\gamma}\leq -L\, C  (r^2+1/n)^{-\gamma/2}\, ,
$$
and then,
by uniform ellipticity, it follows that                  
                  $$ F( D^2 w)- \beta w r^{-\gamma} \leq    (-L\, C  + \beta b^- ) (r^2+1/n)^{-\gamma/2} \leq f(r)(r^2+1/n)^{-\gamma/2}$$
                   as soon as $C$ is chosen large enough. 
  Analogously, for some convenient positive constant $L$,  the function $-L(1-r^\tau)+ b$ is  a radial sub-solution of problem \eqref{dpapp}. 

The standard comparison principle then implies that the sequence $\{ u_n\}$ is uniformly bounded in $C(\overline{B(0,1)})$.  Hence,  it is locally uniformly bounded in $C^2(B(0,1)\setminus\{0\})$ and, up to a subsequence, it is converging locally uniformly in $\overline{B(0,1)}\setminus \{0\}$ to a radial solution $u$ of problem \eqref{dp}, which is a globally bounded function belonging  to $C^2(B(0,1)\setminus \{0\})$.  

Let us now show that the constructed bounded radial solution $u$ is actually continuous in the whole ball $\overline{B(0,1)}$. Indeed, the same approximation argument used in order to prove the existence of $u$ can be applied, in particular,  in order to show the existence of a radial bounded solution $\bar w$  of the Dirichlet problem
$$
\left\{
\begin{array}{cl}
 {\cal M}^+ ( D^2 \bar w) = -(\|f\|_\infty+B+1) r^{-\gamma} & \quad \hbox{ in } B(0,1)\setminus \{0\}\\
 \bar w = u & \quad \hbox{ on } \partial B(0,1)
 \end{array}
 \right.$$
where $B>0$ is a constant such that $\beta |u|\leq B$ in $B(0,1)$.  Lemma \ref{fabiana} and Remark \ref{lemmaformm} applied to $\bar w$ yield that 
                 $$ | \bar w^\prime| \leq c r^{1-\gamma}$$
for some $c>0$.                  
 As a consequence, we have 
                  $$ { \cal M}^+ ( D^2 (u-\bar w) )\geq F( D^2 (u-\bar w) )\geq F(D^2u)-\mathcal{M}^+(D^2\bar w)\geq  r^{-\gamma}\qquad \hbox{ in } B(0,1)\setminus \{0\}$$
 and then,  by Lemma \ref{fabiana} (i),   in a right neighborhood of zero one has
                   $(u-\bar w)^\prime(r) \geq   0\, .$
Hence, 
$$u^\prime(r) \geq -c  r^{1-\gamma}$$
 for $r$ small enough. Analogously, we have
 $$
{ \cal M}^- ( D^2 (u+\bar w) )\leq F( D^2 (u+\bar w) )\leq F(D^2u)+\mathcal{M}^+(D^2\bar w)\leq - r^{-\gamma}\qquad \hbox{ in } B(0,1)\setminus \{0\}$$
 which implies, by Remark \ref{lemmaformm}, 
  $$u^\prime(r) \leq -{\bar w}'\leq  c  r^{1-\gamma}$$
 for $r$ small enough.  Arguing as in the proof of Lemma \ref{fabiana}, from the estimate $|u'(r)|\leq c r^{1-\gamma}$ for $r$ sufficiently small, we deduce that $u$ is  Lipschitz continuous in $\overline{B(0,1)}$ if $\gamma \leq 1$,  it belongs to  $C^1(\overline{B(0,1)})$ if $\gamma <1$, and it is  H\"older continuous in $\overline{B(0,1)}$ with exponent $2-\gamma$ if $\gamma >1$.

Let us observe that the argument above shows that any bounded radial solution of problem \eqref{dp} is continuous in $\overline{B(0,1)}$. This, jointly with
 Theorem \ref{compa}, implies that  problem \eqref{dp} has a unique  radial bounded solution.              
                     
 \end{proof}

The argument used in the above proof yields also the following  compactness result, which we state separately for the sake of clarity.
                   
 \begin{theo}\label{compact}
 Let $\{u_n\}_n$ be a  uniformly bounded sequence of radial functions belonging to $C^2(B(0,1)\setminus \{0\})$ and satisfying
              $$ F( D^2 u_n) = { f_nr^{-\gamma} }\qquad \hbox{ in } B(0,1)\setminus \{0\}\, ,$$ 
               where   $\{f_n\}_n$  are radial, bounded and  continuous on $B(0,1)\setminus \{0\}$.
                If $\{f_n\}$ is uniformly bounded, then $\{u_n\}_n$ are  equicontinuous, thus uniformly converging in $\overline{B(0,1)}$ up to a subsequence. 
                If $\{f_n\}$ is uniformly converging to $f\in C(\overline{B(0,1)})$, then,  up to a subsequence, $\{u_n\}_n$ is  uniformly converging to  a    radial solution $u\in C^2(B(0,1)\setminus \{0\})\cap C(\overline{B(0,1)})$ of 
              $$F( D^2 u) = { r^{-\gamma} f }\qquad  \hbox{ in } B(0,1)\setminus \{0\}\, .$$
              \end{theo}

In the next results, we prove several properties of the \lq\lq smooth" eigenvalue
$$
\bar \lambda_\gamma' \, := \sup \{\mu\, : \  \exists\,   
        u \in C^2(B(0,1)\setminus \{0\})\, ,\ u>0 \hbox{ in } B(0,1)\setminus \{0\}, \ u \hbox{ radial},  \ \ F (D^2 u) + \mu \frac{u}{ r^\gamma} \leq 0\}
 \,   . $$  
In Subsection 2.3 we will prove in facts that $\bar \lambda_\gamma'=\bar \lambda_\gamma$.

Let us start by proving  the validity of the maximum principle  below the value $\bar \lambda_\gamma'$.

\begin{theo}\label{maxpgamma}
 Let    $\mu < \bar \lambda_\gamma'$    and suppose that $u\in C(\overline{B(0,1)})\cap C^2(B(0,1)\setminus \{0\})$ is  a radial function satisfying 
 $$F(D^2 u) + \mu ur^{-\gamma} \geq 0 \qquad \hbox{ in } \ B(0,1)\setminus \{0\}\, .  $$
 If $u(1)\leq 0$,    then $u\leq 0$ in $\overline{B(0,1)}$.
                                   \end{theo}
                     
     \begin{proof}
     If $\mu < 0$, we just apply Theorem \ref{compa} with $v\equiv  g\equiv f \equiv 0$. So, we can assume without loss of generality that $\mu \geq 0$.   
             
 For $\mu^\prime \in ]\mu, \bar \lambda_\gamma'[$,  let $v\in C^2(B(0,1)\setminus \{0\})$ be a radial function satisfying
      $$F( D^2 v) + \mu^\prime  v r^{-\gamma}  \leq 0\, ,\quad v>0\qquad \hbox{ in } B(0,1)\setminus \{0\}\, .$$     We can assume without loss of generality that  $v>0$ on $\partial B(0,1)$, e.g. by performing a dilation in $r$  (this may change a little $\mu^\prime$ but we can still suppose by continuity that $\mu^\prime \in ]\mu, \bar \lambda_\gamma'[$).   
     Arguing as in the proof of Lemma \ref{fabiana}, since ${ \cal M}^-(D^2 v) + \mu' v r^{-\gamma} \leq 0$,  we easily obtain that   $v^\prime(r)$  has constant sign near zero. Assuming by contradiction that $v^\prime (r)\geq 0$ in a neighborhood of zero,  then $v$ is bounded in $B(0,1)\setminus \{0\}$. Hence,  Remark \ref{lemmaformm} applies and  yields $v'(r)\leq 0$ for $r$ small enough: a contradiction. This shows  that $v^\prime(r) \leq 0$ in a neighborhood of $0$. 
     
 Then, there are two possible cases: either $\lim_{r\to 0}v(r)=+\infty$ or $v$ can be extended as a continuous function on $\overline{B(0,1)}$. In the first case, by applying the standard comparison principle, it is easy to prove that for all $\epsilon>0$ one has $u\leq \epsilon v$ in $\overline{B(0,1)}\setminus\{0\}$. In this case, letting $\epsilon\to 0$, we get the conclusion. On the other hand, if $v$ is bounded and continuous on $\overline{B(0,1)}$, we can argue by contradiction. Let us assume that $u$ is positive somewhere in $B(0,1)$, so that   $\frac{u}{v}$  has a positive maximum on $\overline{B(0,1)}$, achieved at some point inside $B(0,1)$. Up to a multiplicative constant for $v$, 
  we  can suppose that 
  $$\max_{\overline{B(0,1)}} \frac{u}{v}=1\, ,$$  
so that  $ u(r)\leq v(r)$. If the maximum is achieved at $0$, then one has 
       $u(0)= v(0)>0\, . $
  By continuity, for   $r$  small enough one has 
       $${\cal M}^+( D^2(u-v) )\geq F(D^2u)- F(D^2v)\geq \frac{1}{2} ( \mu^\prime v(0)-\mu u(0)) r^{-\gamma}\, .$$
Since $\mu^\prime v(0)-\mu u(0)= (\mu'-\mu)v(0) >0$, we can  use Lemma \ref{fabiana} (ii) and  we get that  
    $(u -v)' (r) >0$ in a right neighborhood of $0$. This is a contradiction to  $u-v$ has a maximum point at zero.
 
 Hence, we have that   
  $1= \max \frac{u}{v} > \frac{u(0)}{v(0)}$. Let us select $\eta  <1$ such that 
     $\eta  >\max\{  \frac{\mu}{\mu^\prime}, \frac{u(0)}{v(0)}\}$.  Then,  the function $u-\eta v$ has a  positive maximum achieved at some point  $0<\bar r <1$.     
     
 Since $D^2u(\bar r)\leq \eta D^2v(\bar r)$, by ellipticity  and using \eqref{posh}, we get
      $$ -\mu  v( \bar r){\bar r}^{-\gamma}\leq -\mu  u( \bar r) {\bar r}^{-\gamma} \leq F( D^2 u(\bar r))\leq F( \eta D^2  v (\bar r)) \leq  -\mu^\prime  \eta \, v( \bar r) {\bar  r}^{-\gamma}\, ,$$
which gives the  contradiction $\mu\geq \mu' \eta$. 
       
        \end{proof}
            
The next result provides the existence, uniqueness  and regularity of solutions   below the eigenvalue $\bar \lambda_\gamma'$.

\begin{theo}\label{exilambda}
 Let   $f$ be a  radial and continuous function in $\overline{B(0,1)}$ satisfying $f\leq 0$, with $f$ not identically zero.   Then, for every $\mu < \bar \lambda_\gamma'$ there exists a unique, bounded, radial function $u\in C^2(B(0,1)\setminus \{0\})$ satisfying
 $$
 \left\{
 \begin{array}{cl} 
  F( D^2 u)+ \mu u  r^{-\gamma}  =  f (r) r^{-\gamma} & \hbox{ in } B(0,1)\setminus \{0\}\\
   u=0 & \hbox{ on } \partial B(0,1)
 \end{array}\right.  $$
Moreover, $u$ can be extended as a strictly positive continuous function in $B(0,1)$,   Lipschitz continuous in $\overline{B(0,1)}$ if $\gamma \leq 1$, $(2-\gamma)-$H\"older continuous in $\overline{B(0,1)}$ if $\gamma >1$. 
                \end{theo}
 \begin{proof} 
 
 As in the proof of Theorem \ref{maxpgamma}, we can assume without loss of generality that $\mu>0$, otherwise the conclusion just follows from Theorem \ref{exi1}. The uniqueness of $u$ follows from Theorem \ref{maxpgamma}.
 
  As far as existence is concerned,  let us recursively define a sequence $\{u_n\}_{n\geq 0}$ as follows: we set
$$u_0 \equiv  0\, ,$$
and then, by using Theorem \ref{exi1}, we  define $u_{n+1}\in  C^2(B(0,1)\setminus \{0\}) \cap C(\overline{B(0,1)})$ as the unique bounded radial solution of
$$
\left\{
\begin{array}{cc}
 F( D^2 u_{n+1})  = ( f -\mu u_n )  r^{-\gamma}  & \hbox{ in } B(0,1)\setminus \{0\}\\
 u_{n+1}=0 & \hbox{ on }\partial B(0,1)
 \end{array}\right.$$
  By Theorem \ref{compa}, we have that  $u_{n+1} \geq 0$, hence it is strictly positive in $B(0,1)\setminus \{0\}$  by the standard strong maximum principle, since $f$ is not identically zero. In particular,  $u_n$ is not identically zero for all $n\geq 1$. By applying  the comparison principle in Theorem \ref{compa} again,  we deduce also that   $u_{n+1} \geq u_n$. Let us prove that $\{u_n\}_n$ is uniformly bounded. If not, by setting
                $v_n = {\|u_n\|_\infty}^{-1} u_n $ and $k_n =  \|u_{n+1}\|_\infty^{-1}\|u_n\|_\infty\leq 1$, one gets that $v_{n+1}$ satisfies
$$
\left\{
\begin{array}{cc}
 F ( D^2 v_{n+1}) = \left( \frac{f(r)}{\|u_{n+1}\|_\infty}  -\mu k_n v_n(r)\right)  r^{-\gamma}& \hbox{ in } B(0,1)\setminus \{0\}\\
 v_{n+1}=0 & \hbox{ on }\partial B(0,1)
 \end{array}\right.
  $$
                  Since $\{v_n\}_n$ is uniformly bounded, by applying Theorem \ref{compact},  we can extract  a subsequence still denoted by $\{v_n\}_n$ uniformly converging to  a function $v\geq 0$ satisfying
 $$
 \left\{
\begin{array}{cc}
F ( D^2 v)+ \mu k v r^{-\gamma} =0 & \hbox{ in } B(0,1)\setminus \{0\}\\
 v=0 & \hbox{ on }\partial B(0,1)
 \end{array}\right.
$$
 where $k \leq 1$ is the limit  of some converging subsequence of  $\{k_n\}_n$.
Since $v$ is a radial solution, one has that  $v\in C^2 ( B(0,1)\setminus\{0\})$  and, since $\mu k \leq \mu < \bar \lambda_\gamma'$,  Theorem \ref{maxpgamma} yields $v\leq 0$. Hence, we get $v\equiv 0$, a contradiction with $\|v\|_\infty = 1$. 
                 
                  We have obtained that $\{u_n\}_n$ is bounded, and using once more  Theorem \ref{compact}, we deduce that $\{u_n\}$ uniformly converges to some $u$, which  satisfies the desired equation. 
            By the strong maximum principle, we get that $u>0$ in $B(0,1)\setminus\{0\}$. Moreover, by Remark \ref{lemmaformm}, we have that     $u^\prime (r)\leq 0$ for $r>0$ small enough, which implies  $u(0)>0$.  Finally,    the global regularity of $u$ follows from Theorem \ref{exi1}.
                     
\end{proof} 

We can now prove  that the smooth eigenvalue $\bar \lambda_\gamma'$ is actually achieved on smooth eigenfunctions.
\begin{theo}\label{exieig3}
                     There exists $u\in C(\overline{B(0,1)})\cap C^2(B(0,1)\setminus \{0\})$,   radial,  strictly positive  in $B(0,1)$  and satisfying
$$
 \left\{
\begin{array}{cl}
 F( D^2 u) + \bar \lambda_\gamma' u   r^{-\gamma} = 0& \hbox{ in } B(0,1)\setminus \{0\}\\
 u=0 & \hbox{ on }\partial B(0,1)
 \end{array}\right.
 $$
Furthermore, in $\overline{B(0,1)}$, $u$ is   Lipschitz continuous  when $\gamma \leq 1$ and  H\"older continuous  with exponent $2-\gamma$ if $\gamma >1$.  
\end{theo}
\begin{proof} 
We consider a sequence $\{\lambda_n\}$, with  $\lambda_n \rightarrow \bar \lambda_\gamma'$ and  $\lambda_n <\bar 
\lambda_\gamma'$ and, for all $n$, the solution
 $u_n\in C(\overline{B(0,1)})\cap C^2(B(0,1)\setminus \{0\})$  provided by Theorem \ref{exilambda}
                          of 
                           $$ F( D^2 u_n) + \lambda_n u_n  r^{-\gamma} = {- r^{-\gamma} },\quad  u_n(1)=0 \ .$$
We claim that the positive sequence $\{ \|u_n\|_\infty\}_n$ is unbounded.     Indeed, arguing by contradiction,  if $\{ u_n\}_n$ is uniformly bounded, then,  by using Theorem \ref{compact} and considering a subsequence if necessary,                            we obtain  that there exists a  solution $u\in C(\overline{B(0,1)})\cap C^2(B(0,1)\setminus \{0\})$, $u\geq 0$, of
                            $$ F( D^2 u) + \bar \lambda_\gamma' u  r^{-\gamma}  = {- r^{-\gamma} },\quad  u(1)=0\, . $$
 Then, arguing as in the proof of Theorem \ref{exilambda}, we deduce that $u$ is strictly positive in $B(0,1)$, and by taking $0< \epsilon 
 < \frac{1}{\|u\|_\infty}$,  we see that $u$ satisfies  
                            $$ F( D^2 u) + (\bar \lambda_\gamma'  +\epsilon)  u   r^{-\gamma}     \leq 0, $$
                             a contradiction to the definition of $\bar \lambda_\gamma'$.  It then follows that the sequence $\{u_n\}_n$ is not uniformly bounded.  Normalizing and considering  a subsequence, letting $n\to \infty$ yields the existence of  $u\in  
 C(\overline{B(0,1)})\cap C^2(B(0,1)\setminus \{0\})$ satisfying 
$$ F( D^2 u) + \bar \lambda_\gamma' u  r^{-\gamma} = 0,  \ {\rm in} \   B(0,1) \setminus \{0\} \      , u(1)=0 \ .$$
Finally,  the strict positivity of $u$ in $B(0,1)$ and its global regularity in $\overline{B(0,1)}$ follow by arguing as in the proof of Theorem \ref{exilambda}.

                             \end{proof} 

The uniqueness, up to positive multiplicative constants, of the solution given by Theorem \ref{exieig3} is provided by the last result of the present section.

\begin{prop}\label{simple}
The eigenvalue $\bar \lambda_\gamma'$ is simple. 
\end{prop}
\begin{proof} 

From Theorem \ref{exieig3}, there exists a bounded eigenfunction $v$. Let $u$ be another eigenfunction. 
 
We define 
   $\eta \, : = \sup \{ t>0, t v <u\}$, or, equivalently,  $\frac{1}{\eta} = \sup_{B(0,1)\setminus \{0\}} \frac{v}{u}$, which is well defined by Hopf Lemma applied to $u$ and $v$ on $\partial B(0,1)$. Then, the function $u-\eta v$ is nonnegative in $\overline{B(0,1)}\setminus \{0\}$ and it satisfies
 $$ {\cal M}^-( D^2 ( u-\eta v)) \leq 0\qquad \hbox{ in } B(0,1)\setminus \{0\}\, .
 $$
 If, by contradiction,  $u-\eta v>0$ at some point in $B(0,1)\setminus \{0\}$, then, by the strong maximum principle, $u-\eta v>0$ in the whole of 
 $B(0,1)\setminus \{0\}$. 
 
 We now distinguish the cases  $u(0)$ finite or infinite.

-  In the case $u(0)=+\infty$,  one  necessarily has 
 $$\frac{1}{\eta}=\lim_{r\rightarrow 1} \frac{v}{u} = \frac{v^\prime (1)}{u^\prime (1)}$$
  but this  contradicts Hopf's  Lemma.

- If $u(0) < \infty$, we have either  $\frac{1}{\eta}=\lim_{r\rightarrow 1} \frac{v}{u} = \frac{v^\prime (1)}{u^\prime (1)}$  or $\frac{1}{\eta}=\frac{v(0)}{u(0)}$.  
  Here, the  contradiction follows either from  Hopf's Lemma or from Remark \ref{lemmaformm}, which gives $(u-\eta v)'(r)\leq 0$ for $r>0$ sufficiently small, so that $u-\eta v$ cannot have a strict minimum point at zero.

 Thus,  we have obtained that  all the  eigenfunctions are  bounded and multiple of each others. 
  
\end{proof}

  \subsection{The eigenvalue inherited from some equation on $\R^+$}
  We suppose in this section that $F$ is one of  Pucci's operators.  We consider only the case  $F=\mathcal{M}^+$, the changes to bring for $F=\mathcal{M}^-$ being obvious. 
  
We present below  an alternative proof  of the existence of  radial eigenfunctions related to the eigenvalue $\bar \lambda_\gamma$. Here, the idea is to show   the existence of global solutions defined in  $(0,+\infty)$ of the ODE associated with  radial solutions of  the equation  
       $$ {\cal M}^+ ( D^2 u) = - \frac{u}{r^\gamma}\, \qquad \hbox{ in } \R^N\setminus \{0\}\, .$$  
We recall that $u\in C^2(\R^N\setminus \{0\})$ is a radial solution of the above equation if and only if $u=u(r)$ is a $C^2((0,+\infty))$ solution of the second order ODE
\begin{equation}\label{ODE}
u''=  M_+\left( - \frac{(N-1)}{r} K_+ (u') -\frac{u}{r^\gamma}\right) \qquad \hbox{ in } (0,+\infty)\, ,
\end{equation}
where
$$
K_+(s)=\left\{ \begin{array}{ll}
\Lambda\, s & \hbox{ if } s\geq 0\\
\lambda\, s & \hbox{ if } s< 0
\end{array} \right.\, ,
\qquad 
M_+(s)=\left\{ \begin{array}{ll}
\frac{1}{\Lambda}\, s & \hbox{ if } s\geq 0\\[1ex]
\frac{1}{\lambda}\, s & \hbox{ if } s< 0
\end{array} \right.\, .
$$
        \begin{theo}\label{R}
 There exists  a  global solution $u\in C^2(0,+\infty))$ of  equation \eqref{ODE}, which extends 
        as a continuous function on $[0,+\infty)$ satisfying $u(0)=1$. Moreover, there  exists  $\bar r>0$ such that $u( \bar r) = 0$ and $u(r)>0$ for $0\leq r< \bar r$. 
    \end{theo}
     \begin{proof}
              We begin by proving the local existence of $u$ near zero. We distinguish the cases $\gamma <1$ and $\gamma \geq 1$.        
 \smallskip
        
{\bf $1^{\hbox{st}}$ case} :  $\gamma <1$
         
For fixed $r_0>0$ to be conveniently chosen, let us define the function set   
          $$V_{r_0} = \{ u \in C ( [0, r_0])\, : \  |u(r)-1| \leq \frac{1}{2}, u(0) = 1\}. $$  
For $u\in V_{r_0}$, let us define                    
 $$ T(u)(r)\, :  = 1-\int_0^r \frac{1}{\lambda s^{N-1}} \int_0^s u(t) t^{N-1-\gamma} dt ds.$$
We fix $r_0$ such that $r_0^{2-\gamma}\leq\frac{\lambda (N-\gamma)(2-\gamma)}{3}$. Then, it is easy to verify that   $T$  maps $V_{r_0}$ into itself and it is  a contraction mapping on it. Let us denote by $u\in V_{r_0}$  its fixed point. 
            
  It then follows that $u$ satisfies
   $$u^\prime (r) = -\frac{1}{\lambda  r^{N-1}} \int_0^r u(t) t^{N-1-\gamma} dt$$
   as well as
   $$
u''=-\frac{u}{\lambda r^\gamma} -( N-1) \frac{u^\prime}{r}\, .
$$
Thus, we clearly have $u'\leq 0$ and,  if we prove that $u''\leq 0$ as well, then $u$ is a solution of \eqref{ODE} in $(0,r_0)$.
We observe that   
  $u^\prime \leq 0$ implies $u\leq 1$ and,  consequently, 
   $$ u^\prime (r) \geq -\frac{1}{\lambda r^{N-1}} \int_0^r t^{N-1-\gamma} dt = -\frac{r^{1-\gamma}}{\lambda ( N-\gamma)}\, .$$  
Hence 
\begin{equation}\label{ubelow}
 u(r) \geq 1-\frac{r^{2-\gamma}}{\lambda(2-\gamma)  ( N-\gamma)}\, ,
\end{equation} 
which  in turn  implies 
     \begin{eqnarray*}
      u^\prime (r) &\leq& -\frac{1}{\lambda r^{N-1}} \int_0^r (1-\frac{t^{2-\gamma}}{\lambda(2-\gamma)  ( N-\gamma)})t^{N-1-\gamma} dt\\
      &=&  -\frac{r^{1-\gamma}}{\lambda ( N-\gamma)}+ \frac{r^{3-2\gamma}}{\lambda^2(2-\gamma)  ( N-\gamma)(N+2-2\gamma)}\\
      & \leq&  -\frac{r^{1-\gamma}}{2\lambda ( N-\gamma)}\end{eqnarray*}
       by the choice of $r_0$.  Thus, we have proved that
             \begin{equation}\label{encaduprime}
              -\frac{r^{1-\gamma}}{\lambda ( N-\gamma)} \leq u^\prime (r)
\leq     -\frac{r^{1-\gamma}}{2\lambda ( N-\gamma)} 
\end{equation}
From estimates \eqref{ubelow} and \eqref{encaduprime}, we further deduce
 \begin{eqnarray*}
       u^{\prime \prime } &=& -\frac{u}{\lambda r^\gamma} -(N-1) \frac{u^\prime}{r} \\
       &\leq& -\left( 1-\frac{r^{2-\gamma}}{\lambda(2-\gamma) ( N-\gamma)}\right)\frac{r^{-\gamma}}{\lambda}  + (N-1) \frac{r^{-\gamma}}{\lambda ( N-\gamma)} \\
       &= & -\frac{(1-\gamma)r^{-\gamma}}{\lambda (N-\gamma)}  + \frac{r^{2-2\gamma}}{ \lambda^2(2-\gamma)  ( N-\gamma)}\\
       &\leq & -\frac{(1-\gamma)r^{-\gamma}}{2\lambda(N-\gamma)} 
       \end{eqnarray*}
 if we fix $r_0$ by setting $r_0^{2-\gamma}=\frac{\lambda (1-\gamma)(2-\gamma)}{3}$. By this choice for $r_0$, we obtain that  $u$ satisfies equation \eqref{ODE}
  in $(0, r_0]$ and, moreover, that $u\in C^1([0,r_0])$.
  \smallskip
  
 {\bf $2^{\hbox{nd}}$ case} :  $\gamma \geq1$         
 
 In this case, the solution $u$ is expected to be locally convex near zero. Thus, for $v\in V_{r_0}$, we define the map
               $$ T(v)(r) = 1-\int_0^r \frac{1}{\Lambda s^{\tilde{N}_+-1}} \int_
             0^s u(t) t^{{\tilde{N}_+-1}-\gamma} dt\,  ds$$
As in the previous case, it is easy to check that if $r_0^{2-\gamma}\leq \frac{\Lambda (\tilde{N}_+-\gamma)(2-\gamma)}{3}$,  
 then $T$ is a contraction mapping on $V_{r_0}$.            
             Let $u\in V_{r_0}$ be its fixed point. We then have
$$u^\prime (r) = -\frac{1}{\Lambda r^{ \tilde N_+-1}} \int_
             0^r u(t) t^{ \tilde N_+-1-\gamma} dt $$   
and
$$
u''=-\frac{(\tilde{N}_+-1)}{r}  u'-\frac{u }{\Lambda r^{\gamma}}\, .
$$
Thus, in order to prove that $u$ satisfies \eqref{ODE} in $(0,r_0]$, it is enough to show that $u''\geq 0$.
                     
By using a bootstrap argument analogous to the one used in the first case, we deduce             
$$ 
-\frac{r^{1-\gamma}}{\Lambda (\tilde{N}_+-\gamma)} \leq u^\prime \leq   -\frac{r^{1-\gamma}}{\Lambda (\tilde{N}_+-\gamma)} +\frac{r^{3-2\gamma}}{\Lambda^2 (\tilde{N}_+ -\gamma)(2-\gamma)(\tilde{N}_++2-2\gamma)}
 $$
as well as
$$
 1- \frac{r^{2-\gamma}}{\Lambda (\tilde{N}_+-\gamma)(2-\gamma)}\leq u \leq 1- \frac{r^{2-\gamma}}{\Lambda (\tilde{N}_+-\gamma)(2-\gamma)} +\frac{r^{4-2\gamma}}{2\Lambda^2 (\tilde{N}_+ -\gamma)(2-\gamma)^2(\tilde{N}_++2-2\gamma)}\, .
$$
We then have
                 \begin{eqnarray*}
                  u^{\prime \prime }(r) &=  & -\frac{u}{\Lambda r^\gamma} -\frac{(\tilde{N}_+-1)}{r}u' \\
&\geq &  \frac{\gamma-1}{\Lambda (\tilde{N}_+-\gamma)} r^{-\gamma}
+\frac{3-2\gamma}{\Lambda^2 (\tilde{N}_+ -\gamma)(2-\gamma)(\tilde{N}_++2-2\gamma)} r^{2-2\gamma}\\
& & 
-\frac{r^{4-3\gamma}}{2\Lambda^3 (\tilde{N}_+ -\gamma)(2-\gamma)^2(\tilde{N}_++2-2\gamma)}
  \end{eqnarray*}
By choosing $r_0$ such that
$$
r_0^{2-\gamma}\leq 
\Lambda (3-2\gamma)(2-\gamma) 
$$
we obtain in any case that $u''\geq 0$ in $(0,r_0]$.
\smallskip
             
Thus, there exists a local solution $u$ of equation \eqref{ODE} in $(0,r_0]$, which is positive and satisfies the initial condition   $u(0)=1$. By  observing the Lipschitz continuity of the functions $K_+$ and $M_+$ and using Cauchy--Lipschitz Theorem, we can extend the solution $u$ as a global solution in $(0,+\infty)$.                                   
             
It remains               to prove that there exists a first point $\bar r$ such that  $u( \bar r) = 0$ and $u(r)>0$ for $0\leq r<\bar r$. 

We argue by contradiction, and we assume that $u(r)>0$ for all $r\geq 0$. Since $u'$ is initially negative, if there exists a first point $r_1>0$   such that $u'(r_1)=0$, then $u''(r_1)\geq 0$. On the other hand, from equation \eqref{ODE}, we deduce $u''(r_1)<0$ since $u(r_1)>0$. Hence, one has $u'(r)<0$ for all $r>0$. From \eqref{ODE}, it then follows that, independently of the sign of $u^{\prime \prime}$,   one has
         $$ u^{\prime \prime} +  \frac{N-1}{r} u^\prime \leq -\frac{u r^{-\gamma}}{\Lambda}\, .$$
         Inspired by \cite{BEQ} and \cite{D}, let us introduce the function
           $$ y(r) = \frac{u^\prime(r)}{u(r)} r^{N-1}\, ,$$
which  is, then,  negative and it satisfies 
            $$ y^\prime \leq  -\frac{r^{N-1-\gamma}}{\Lambda}-\frac{y^2}{r^{N-1}}\, .$$
By integrating  between some $r_1>0$ and $r$ , we obtain 
            $$ y(r) + k(r) \leq -c_1 r^{N-\gamma}\, ,$$
for some $c_1>0$ and $k(r)=\int_{r_1}^r \frac{y^2(t)}{t^{N-1}}\, dt$.
    This yields   in particular   $y(r)\leq -c_1 r^{N-\gamma}$ and therefore, for $r$ sufficiently large,
\begin{equation}\label{k1}
k(r)\geq c_2 r^{N+2(1-\gamma)}\, .
\end{equation}
 On the other hand, we also have 
              $$k(r) \leq -y(r)\, ,$$
              that is
 $$k(r)\leq \sqrt{ k'(r)r^{N-1}}\, $$
 which yields, after integration on $(r,+\infty)$, 
\begin{equation}\label{k2}
k(r)\leq (N-2) r^{N-2}\, .
\end{equation}
Being $N+2(1-\gamma)>N-2$, estimates \eqref{k1} and \eqref{k2} give a contradiction, showing that the constructed solution $u$ cannot be globally  positive in $[0,+\infty)$.
              
 \end{proof}  
 
As  an immediate consequence of the above result and  Theorem \ref{maxpgamma}, we deduce the following
 \begin{cor}\label{barlambda}
Let $\bar r$ be defined as in Theorem \ref{R}. Then   $\bar \lambda_\gamma' ( B(0,1)\setminus\{0\} ) = \bar r^{\gamma-2}$. 
\end{cor}

                        \begin{rema}
{\rm                        Let us observe that, as in the case of equations with  continuous   coefficients, one could prove the existence of a numerable set of  radial eigenvalues, by  proving the  oscillatory behavior of the solution $u$ constructed in Theorem \ref{R}, 
                         see \cite{BEQ} and \cite{D}.}
                        \end{rema}   
                               
                                \subsection{The stability of the principal eigenvalue and  related eigenfunctions} 
 The results of the present  section give, as a corollary,  the proof of Theorem \ref{exigamma}. 
 
  Let us start by proving the stability with respect to the $\epsilon-$regularization of the singular potential. We recall that $r_\epsilon = (r^2+\epsilon^2)^\frac{1}{2}$ and $\bar \lambda_\gamma^\epsilon = \bar \lambda ( F, \frac{1}{r_\epsilon^\gamma}, B(0,1))$. 
 
  \begin{theo}\label{th3}
  One has
       $$\bar \lambda_\gamma' =\lim_{ \epsilon \rightarrow 0} \bar \lambda_\gamma^\epsilon\, .$$ 
        Furthermore,  if $\{u_\epsilon\}$ is  the sequence of the eigenfunctions associated with the eigenvalue $\bar \lambda_\gamma^\epsilon$ and satisfying  $u_\epsilon (0) = 1$, then, one can extract from $\{u_\epsilon\}$ a subsequence uniformly converging on $\overline{ B(0,1)}$ to the eigenfunction associated with $\bar \lambda_\gamma'$ which takes the value $1$ at zero.
 \end{theo}
        
           \begin{proof} 
 Let   $\{u_\epsilon\}$ be  the sequence as in the statement. Then, each $u_\epsilon$ is a smooth positive function in $B(0,1)$ satisfying in particular
 $$
F(D^2u_\epsilon)+\bar \lambda_\gamma^\epsilon \frac{u_\epsilon}{r^\gamma}\geq 0 \quad \hbox{ in } B(0,1)\setminus \{0\}\, ,
 $$
 so that, by Theorem \ref{maxpgamma}, one has $\bar \lambda_\gamma^\epsilon\geq \bar \lambda_\gamma'$. Moreover, the sequence $\{ \bar \lambda_\gamma^\epsilon\}$ is monotone increasing with respect to $\epsilon$. Thus, we deduce
 $$
 \mu\, : = \lim_{\epsilon\to 0} \bar \lambda_\gamma^\epsilon \geq \bar \lambda_\gamma'\, .$$
On the other hand, by the monotonicity properties of radially symmetric solutions of elliptic equations, we know that $u_\epsilon'(r)\leq 0$ for $r\in [0,1]$. 
Since $${ \cal M}^+ ( D^2 u_\epsilon) \geq F( D^2 u_\epsilon)$$
we deduce that, independently of the sign of $u_\epsilon''(r)$, one has
$$ u_\epsilon^{\prime \prime} + (\tilde{N}_+-1) \frac{u_\epsilon^\prime}{r} \geq -\frac{\bar \lambda_\gamma^\epsilon}{\lambda} u_\epsilon  r_\epsilon^{-\gamma} \, .$$
This implies
$$ ( u_\epsilon^\prime r^{\tilde{N}_+-1})^\prime \geq - \frac{\bar \lambda_\gamma^\epsilon}{\lambda} u_\epsilon \frac{r^{\tilde{N}_+-1}} {r_\epsilon^\gamma}\geq - \frac{\bar \lambda_\gamma^\epsilon}{\lambda} r^{\tilde{N}_+-1-\gamma}\, ,$$ 
and therefore, 
by integrating, 
$$
 0\geq u_\epsilon^\prime(r)\geq - \frac{\bar \lambda_\gamma^\epsilon}{\lambda (\tilde{N}_+-\gamma)} r^{1-\gamma}\, .
 $$
Hence, on $\overline{B(0,1)}$, the functions $u_\epsilon$ are uniformly Lipschitz continuous if $\gamma \leq 1$, and uniformly $(2-\gamma)-$H\"older continuous if $\gamma>1$. In both cases, up to a subsequence, $\{u_\epsilon\}$ is uniformly converging to a continuous radial function $u\in C(\overline{B(0,1)})$ which satisfies $u(0)=1$ and 
$$ F( D^2 u) + \mu \frac{u}{r^\gamma}=0.$$
Hence,  $u$ is $C^2(B(0,1)\setminus \{0\})$ and, by the standard strong maximum principle, $u$ is strictly positive in $B(0,1)$. This yields, by definition, $\mu\leq \bar \lambda_\gamma^\prime$. Hence, $\mu =\bar \lambda_\gamma^\prime$ and the conclusion follows from Proposition \ref{simple}.

\end{proof}

As a consequence of the previous theorem, we finally obtain the following 

        \begin{cor}\label{delta} One has
         $$\bar\lambda_\gamma = \lim_{\delta\to 0} \bar \lambda _\gamma \left( B(0,1) \setminus \overline{B(0, \delta)}\right) =\bar \lambda_\gamma'\, .$$
          \end{cor}

      \begin{proof} 
 
 We observe that the function     $\delta \mapsto  \bar \lambda_\gamma \left( B(0,1)\setminus B(0, \delta)\right) $ is monotone increasing. Moreover, by their own definition, we have that
 $$
\bar \lambda_\gamma'\leq  \bar \lambda_\gamma \leq \bar \lambda_\gamma \left( B(0,1)\setminus B(0, \delta)\right)\, \quad \hbox{ for all } \delta\geq0\, .
 $$
On the other hand,     by Theorem \ref{th3},  for any $\eta >0$ there exists $\epsilon_0>0$ such that  
   $$ \bar \lambda_\gamma^{ \epsilon_0}  \leq \bar \lambda _\gamma' + \frac{\eta}{2}\, .$$
Furthermore,   by using the continuity of the principal eigenvalue with respect to the domain for equations with regular coefficients, there exists $\delta_0>0$ such that 
   $$\bar \lambda_\gamma^{\epsilon_0} ( B(0,1))\setminus B(0, \delta_0)) \leq \bar \lambda _\gamma ^{\epsilon_0}   + \frac{\eta}{2} \leq  \bar \lambda_\gamma' +  \eta \, .$$
Now,  since  $\epsilon \mapsto \bar  \lambda^\epsilon_\gamma \left(B(0,1))\setminus B(0, \delta_0)\right) $ decreases when $\epsilon$ decreases to zero, one gets 
  $$  \bar \lambda_\gamma \left( B(0,1)\setminus B(0, \delta_0)\right)\leq \bar \lambda_\gamma^{\epsilon_0} \left( B(0,1))\setminus B(0, \delta_0)\right) \leq \bar \lambda_\gamma' + \eta\, ,$$
   which gives the conclusion.
   
 \end{proof} 
   
 \section{ The case $\gamma=2$ : Proof of Theorem \ref{gamma=2}.}
 
 The aim of this section is to give the proof of Theorem \ref{gamma=2}. 
  
 Let us start by recalling that, in the semilinear case, the eigenvalue related to Laplace operator with an inverse quadratic potential can be defined by a variational approach, i.e. by considering the minimum problem
 $$ \bar \lambda_2 ( \Delta)\, :\, = \inf_{ \stackrel{u\in H_0^1 ( B(0,1))} {\int_{B(0,1)} \frac{u^2}{|x|^2} dx = 1}} \int_{B(0,1)} |\nabla u |^2 dx.$$
In this case, one has
$$
\bar \lambda_2 ( \Delta) = \left( \frac{N-2}{2}\right)^2\, .
$$
Indeed, on the one hand, by Hardy inequality, every function $u\in H^1_0(B(0,1)$ satisfies
$$
\int_{B(0,1)} \frac{u^2}{|x|^2} dx \leq \left( \frac{2}{N-2}\right)^2 \int_{B(0,1)} | \nabla u |^2 dx$$
      and then 
     $$ \bar \lambda_2 ( \Delta) \geq \left( \frac{N-2}{2}\right)^2\, .$$
On the other hand, for every $\epsilon>0$, the function    $u_\epsilon (r) = r^{-\frac{N-2}{2} + \epsilon} ( -\log r)$
      belongs to $H_0^1(B(0,1))$ and 
        satisfies 
        $$ \int_{B(0,1)} | \nabla u_\epsilon|^2 = \left[\left(\frac{N-2}{2}\right)^2+  \epsilon^2\right] \int_{B(0,1)} 
        \frac{|u_\epsilon|^2}{|x|^2} dx\, , $$
so that
 $$ \bar \lambda_2 ( \Delta) \leq \left( \frac{N-2}{2}\right)^2+\epsilon^2\, ,\qquad \forall\, \epsilon>0\, .$$
 As it is well known, the infimum defining $\bar\lambda_2(\Delta)$ is not achieved, that is the Dirichlet problem
 $$
 \left\{
 \begin{array}{cl}
 -\Delta u = \bar\lambda_2(\Delta) \frac{u}{r^2} & \hbox{ in } B(0,1)\\[2ex]
 u=0 & \hbox{ on } \partial B(0,1)
 \end{array}
 \right.
 $$
 has not finite energy solutions $u\in H^1_0(B(0,1))$.
 
 A kind of variational approach is possible also in the fully nonlinear framework for the case of Pucci's operators. From now on,  we consider the operator $\mathcal{M}^+$, being obvious the changes to be made for the operator $\mathcal{M}^-$.
 
  Let us introduce the space of functions 
 $$\mathcal{V}= \left\{ u\in C^2([0,1])\, : u'(0)=0\, ,\ {\rm supp}(u) \hbox{ compact in } [0,1)\right\}\, ,
 $$
  endowed with the norm
 $$
 \| u\| =\left( \int_0^1 |u'|^2 r^{\tilde{N}_+-1} dr\right)^{1/2}\, ,
 $$
 and let us denote by $\mathcal{H}^1_0$  the closure of $\mathcal{V}$. Then, for all $\gamma\leq 2$, we can consider the minimum problem
$$
 \bar \lambda_{\gamma,var} \, :\, =  \inf_{ \stackrel{u\in \mathcal{H}_0^1} {\int_0^1 u^2r^{\tilde{N}_+-1-\gamma}dr = 1}} \int_0^1 |u'|^2 r^{\tilde{N}_+-1}dr\, .
 $$ 
 In the next results, we will relate the two values $\bar \lambda_\gamma(\mathcal{M}^+)$ and $\bar \lambda_{\gamma,var}$, and we will study their asymptotic behavior as $\gamma\to 2$.

 \begin{theo}\label{barbar}
One has
       $$\bar \lambda_{2,var}= \left( \frac{\tilde{N}_+-2}{2}\right)^2\, .$$
 \end{theo}
 \begin{proof}     For any $\epsilon>0$,    let  $u_\epsilon(r)= r^{-\frac{\tilde{N}_+-2}{2}+\epsilon } ( -\log r)$. Then, it is easy to check that  $u_\epsilon \in \mathcal{H}_0^1$ and a direct computation shows that
$$ \int_0^1|u_\epsilon ^\prime|^2 r^{\tilde{N}_+-1} dr=  \left[ \left( \frac{\tilde{N}_+-2}{2}\right)^2+\epsilon^2\right] 
\int_0^1  u_\epsilon^2 r^{\tilde{N}_+-3} dr\, ,
$$
hence
$$
\bar \lambda_{2,var}\leq  \left( \frac{\tilde{N}_+-2}{2}\right)^2+\epsilon^2\, , \qquad \forall\, \epsilon>0\, .
$$
On the other hand, we observe that the function $u= r^{-\frac{\tilde{N}_+-2}{2}} ( -\log r)$ satisfies, for $r> 0$ ,
      $$u^{\prime \prime} + (\tilde{N}_+-1)\frac{u'}{r} = - \left( \frac{\tilde{N}_+-2}{2}\right)^2 u r^{-2}.$$ 
   Let us multiply the above equation by $\frac{v^2}{u}r^{\tilde{N}_+-1}$,  where $v \in \mathcal{V}$ is arbitrarily fixed. 
 Since $\tilde{N}_+>2$, we have that $\frac{u^\prime r^{\tilde{N}_+-1}}{u}$ tends to zero  as $r\to 0$. As a consequence,   integrating by parts, we get
 $$-\left( \frac{\tilde{N}_+-2}{2}\right)^2 \int_0^1 v^2 r^{\tilde{N}_+-3}dr =  \int_0^1 \left( \frac{u^\prime}{u} v -v^\prime\right)^2 r^{\tilde{N}_+-1}dr - \int_0^1 ( v^\prime )^2 r^{\tilde{N}_+-1}dr \, ,$$
        which yields, by the arbitrariness of $v\in \mathcal{V}$, 
                $$ \bar\lambda_{2,var}\geq \left( \frac{\tilde{N}_+-2}{2}\right)^2\, .$$
          
      \end{proof}
      
In order to establish the relationship between  $\bar \lambda_{\gamma,var}$ and $\bar \lambda_\gamma$ we need to investigate on the monotonicity and convexity properties of the functions $u$ realizing the infimum in the definition of $\bar \lambda_{\gamma,var}$.

      \begin{prop}\label{signederiveeseconde}
      Let $1<\gamma <2$ and assume that $u_\gamma\in \mathcal{H}^1_0$, with $u_\gamma\geq0$, realizes the infimum defining $\bar \lambda_{\gamma, var}$. Then, $u_\gamma \in C^2 ((0,1])$ is bounded,  $u_\gamma^\prime \leq 0$ and  $u_\gamma^{\prime \prime } \geq 0$ in $(0,1)$.       \end{prop}
      
 \begin{proof} 
      
       Since $u_\gamma$ is a minimum,  for any $v \in \mathcal{H}^1_0$ one has
       $$ \int _0^1u_\gamma ^\prime v^\prime   r^{\tilde{N}_+-1} =  \bar \lambda_{\gamma, var} \int _0^1u_\gamma  v  r^{\tilde{N}_+-1-\gamma}.$$
         In particular, $u_\gamma$ satisfies  in the distributional sense 
         \begin{equation}\label{euler}-(u^\prime_\gamma   r^{\tilde{N}_+-1})^\prime =  \bar \lambda_{\gamma, var}  u_\gamma  r^{\tilde{N}_+-1-\gamma}
         \end{equation} 
By regularity theory, this implies that $u_\gamma$ belongs to $C^2((0,1])$, it is strictly positive in $(0,1)$  and it satisfies $u_\gamma(1)=0$. Let us prove that $u_\gamma$ is bounded and that it can be extended as a continuous function on $[0,1]$. Indeed, by multiplying                  
          equation \eqref{euler} by a smooth function $v\in C^2([0,1])$, having  compact support in $[0,1)$ and satisfying $v(0)\neq0$, and integrating on $[\epsilon, 1]$ for $\epsilon>0$, one has  
 $$
\bar \lambda_{\gamma, var}  \int_\epsilon^1 u_\gamma(r)  r^{\tilde{N}_+-1-\gamma} v(r)\, dr= 
            \int_\epsilon ^1   u^\prime_\gamma(r)   r^{\tilde{N}_+-1} v^\prime(r)\, dr + u^\prime_\gamma  ( \epsilon)   \epsilon ^{\tilde{N}_+-1} v( \epsilon) \, .
$$
Letting $\epsilon$ go to zero,  we deduce 
            $ \lim _{ \epsilon \rightarrow 0}  u_\gamma ^\prime ( \epsilon)   \epsilon ^{\tilde{N}_+-1} = 0$.
It then follows, again from \eqref{euler},  that  $u^\prime_\gamma (r) \leq 0$ and that there exists some  positive constant $c_0$ such that
             $$u^\prime_\gamma(r)  \geq -c_0\, r^{1-\tilde N_+}\qquad \hbox{ for } r\in (0,1]\, .$$
This implies that
$$u_\gamma (r)=-\int_r^1 u'_\gamma(s)\, ds\leq c_0\int_r^1s^{1-\tilde{N}_+} ds\leq d_0 r^{2-\tilde{N}_+}$$
with $d_0=\frac{c_0}{\tilde{N}_+-2}>0$, which, in turn, yields
$$
u'_\gamma(r)r^{\tilde{N}_+-1}=-\bar \lambda_{\gamma, var} \int_0^ru_\gamma(s) s^{\tilde{N}_+-1-\gamma}ds
\geq -\bar \lambda_{\gamma, var} \, d_0 \int_0^rs^{1-\gamma}ds =-c_1 r^{2-\gamma}\, .
$$
Thus, we have
$$u^\prime_\gamma(r)  \geq -c_1\, r^{1-\tilde N_++2-\gamma}
$$
and, then,
 $$
 u_\gamma(r)\leq   d_1   r^{2-\tilde{N}_++2-\gamma}\, .
 $$
 Iterating the above inequalities, we obtain that for all integers $j\geq0$ such that $2-\tilde{N}_++j(2-\gamma)<0$, there exist positive constants $c_j$ and $d_j$ satisfying
\begin{equation}\label{estj}
 u^\prime_\gamma(r)  \geq -c_j\, r^{1-\tilde N_++j(2-\gamma)}\, ,\quad  u_\gamma(r)\leq   d_j   r^{2-\tilde{N}_++j(2-\gamma)}\, .
 \end{equation}
Now, if there exists $j\in \N$ such that $2-\tilde{N}_++j(2-\gamma)=0$, i.e. if $\frac{\tilde{N}_+-2}{2-\gamma}\in \N$, then, by integrating the estimates obtained  at the $(j-1)$-th step, we obtain
$$
u'_\gamma(r)\geq -c_j \, r^{-1}\, ,\quad u_\gamma(r)\leq d_j\, (-\ln r)\, .
$$
Integrating once more, we finally deduce
$$
u'_\gamma(r)\geq -c_{j+1} (-\ln r) r^{1-\gamma}\Longrightarrow u_\gamma(r)\leq d_{j+1}=\frac{c_{j+1}}{(2-\gamma)^2}\, .
$$
On the other hand, if $\frac{\tilde{N}_+-2}{2-\gamma}$ is not integer, by integrating estimates \eqref{estj} for $j=\left[ \frac{\tilde{N}_+-2}{2-\gamma}\right]$, we obtain
$$
u'_\gamma(r)\geq -c_{j+1} r^{1-\tilde{N}_+ +(j+1)(2-\gamma)} \Longrightarrow u_\gamma(r)\leq d_{j+1}=
\frac{c_{j+1}}{2-\tilde{N}_++(j+1)(2-\gamma)}\, .
$$
This shows  that, in any case, $u_\gamma$ is bounded. 

Let us finally  prove that $u_\gamma^{\prime \prime } \geq 0$. We introduce the function
 $$y_\gamma(r) \, : = (\tilde{N}_+-1) u_\gamma^\prime(r) + { \bar \lambda_{ \gamma, var}  } r^{1-\gamma} u_\gamma(r)\, , $$
 which verifies
       $y_\gamma = -r   u^{\prime \prime}_\gamma$. Hence, 
        we need to prove that $y_\gamma (r)\leq 0$  for $r\in (0,1]$. 
        An easy computation shows that
        $$ y_\gamma^\prime (r)+ (\tilde{N}_+-1) \frac{y_\gamma(r)}{r} =   { \bar \lambda_{\gamma, var} }\left( u^\prime_\gamma(r) r^{1-\gamma} +(1-\gamma) u_\gamma(r)  r^{-\gamma}\right) \leq 0\, ,$$
 so that
         $$ (y_\gamma (r) r^{ \tilde N_+-1})^\prime \leq 0\, .$$
Since  $u_\gamma$ is bounded and  $-c_0  r^{1-\gamma}\leq u^\prime_\gamma(r) \leq 0$, we deduce that 
$y_\gamma(r) r^{\tilde{N}_+-1}\rightarrow 0$ as $r\to 0$. Hence,           $y_\gamma(r) \leq 0$ for $r>0$.      
     
           \end{proof}

           \begin{cor}\label{corvargamma}
           Let $\gamma\in ]1,2[$. Then
           $$\bar \lambda_\gamma (\mathcal{M}^+) = \Lambda \, \bar\lambda_{\gamma, var}\, .$$
                      \end{cor}
           
            \begin{proof}
It is not difficult to prove that the infimum defining $\bar \lambda_{\gamma,var}$ is achieved for $1<\gamma<2$. Thus, there exists $v_\gamma\in \mathcal{H}^1_0$,   which can be assumed to be positive in $[0,1)$, and which satisfies in $(0,1)$
            $$ \Lambda v_\gamma^{\prime \prime } + \lambda \frac{N-1}{r} v_\gamma^\prime = - \Lambda \, \bar\lambda_{\gamma, var}v_\gamma r^{-\gamma}.$$
             By Proposition \ref{signederiveeseconde},  we have that  $v_\gamma$ is bounded, $v_\gamma^\prime \leq 0$ and  $v_\gamma^{\prime \prime} \geq 0$, so that $v_\gamma$ satisfies
              $${\cal M}^+( D^2 v_\gamma) +\Lambda\,  \bar \lambda_{\gamma, var}   v_\gamma r^{-\gamma}= 0\quad \hbox{ in } B(0,1)\setminus \{0\}\, .$$
   By definition,   it then follows  that    $\bar \lambda_\gamma(\mathcal{M}^+)=\bar \lambda_\gamma'(\mathcal{M}^+) \geq \Lambda\, \bar\lambda_{\gamma, var}$. Furthermore, analyzing the boundary condition, we get, by regularity, that  $v_\gamma(1) =0$ in the classical sense.   If, by contradiction,  $\Lambda\, \bar\lambda_{\gamma, var}< \bar \lambda_\gamma'(\mathcal{M}^+)$, then Theorem \ref{maxpgamma} would give $v_\gamma \leq 0$ in $B(0,1)$, a contradiction. 
 
             \end{proof}  
      
      \begin{cor} \label{convgammavar} One has
      $$\lim_{\gamma \rightarrow 2}  \bar \lambda_\gamma (\mathcal{M}^+)=  \Lambda \, \left( \frac{\tilde{N}_+-2}{2}\right)^2\, .$$
      \end{cor} 
      \begin{proof} 
           By Theorem \ref{barbar} and Corollary \ref{corvargamma}, it is sufficient to prove that $ \bar \lambda_{\gamma , var}\rightarrow \bar \lambda_{2,var}$ as $\gamma\to 2$. 
           
         We first observe that, by their own definition,  $ \bar \lambda_{2, var}\leq \bar \lambda_{\gamma, var}$.  
         
On the other hand, for  any  $\epsilon>0$ there exists $v\in \mathcal{H}^1_0$ such  that 
         $$  \int_0^1 |v^\prime |^2  r^{\tilde{N}_+-1} dr \leq ( \bar  \lambda_{2,var}+ \epsilon) 
         \int_0^1 |v |^2  r^{\tilde{N}_+-3} dr\, .$$
Moreover,    there exists $\gamma_0$ sufficiently close  to $2$ in order that, for $\gamma \geq  \gamma_0$ , 
$$ \int_0^1 |v |^2  r^{\tilde{N}_+-1 -\gamma } dr \geq  (1-\epsilon) \int_0^1 |v |^2  r^{\tilde{N}_+-3} dr\, .$$
Thus, one has
$$\int_0^1 |v^\prime |^2  r^{\tilde{N}_+-1} dr \leq ( \bar  \lambda_{2,var}+ \epsilon) (1-\epsilon)^{-1} \int_0^1 |v |^2  r^{\tilde{N}_+-1 -\gamma } dr
$$
 which yields
             $$\bar \lambda_{\gamma, var} \leq  ( \bar \lambda_{ 2,var}+ \epsilon) (1-\epsilon)^{-1}\, . $$
 
                \end{proof}

 We are now ready to prove statement (i) of Theorem \ref{gamma=2}.

  \begin{theo}\label{exp}
 One has
 $$
\bar \lambda_2 (\mathcal{M}^+)=    \Lambda\, \left( \frac{\tilde{N}_+-2}{2}\right)^2
$$
and the function $u(r)= r^{-\frac{\tilde{N}_+-2}{2}} (-\ln r)$ is an explicit solution of
$$\left\{ \begin{array}{cl}
 \mathcal{M}^+( D^2 u) + \bar \lambda_\gamma \frac{u}{r^\gamma} = 0 & \hbox{ in } \ B(0,1)\setminus \{0\}\\[1ex]
       u = 0 & \hbox{ on } \ \partial B(0,1)
       \end{array}\right.$$ \end{theo}
               \begin{proof} 
For any positive constants $c_1$ and $c_2$, let us consider the function
   \begin{equation}\label{fonction}u(r)= r^{-\frac{\tilde{N}_+-2}{2}}(c_1 ( -\ln r)+c_2).
   \end{equation}
   An easy computation, analogous to the one made in the proof of Theorem \ref{barbar}, leads to 
$$
    \mathcal{M}^+ ( D^2 u) +  \Lambda \left( \frac{\tilde{N}_+-2}{2}\right)^2  \frac{u}{ r^2} = 0\quad \hbox{ in } B(0,1)\setminus \{0\}\, .
$$
This gives, by definition,  that 
    $$\bar  \lambda_2(\mathcal{M}^+)\geq  \Lambda \left( \frac{\tilde{N}_+-2}{2}\right)^2\, .$$
  On the other hand, by observing that  
              $ \bar \lambda_\gamma \geq \bar \lambda_2$ for all $\gamma\leq 2$ and by using  Corollary \ref{convgammavar}, 
              we also have 
              $$ \bar \lambda_2 (\mathcal{M}^+)\leq \ \lim_{\gamma \to 2} \bar \lambda_\gamma (\mathcal{M}^+)= \Lambda \left( \frac{\tilde{N}_+-2}{2}\right)^2\, .$$
   
              \end{proof}
            
 As a consequence of Corollary \ref{convgammavar} and Theorem \ref{exp}, we immediately deduce the first stability property of $\bar \lambda_2 (\mathcal{M}^+)$ stated in Theorem \ref{gamma=2}-(ii). The other ones are given by 
 by the following result.
         
              \begin{cor}\label{convdelta} For the operators $F=\mathcal{M}^\pm$, one has 
              $$\lim_{ \delta \rightarrow 0} \bar \lambda_2 ( B(0,1) \setminus  \overline{ B(0, \delta)})=\bar \lambda_2 =\lim_{\epsilon\to 0} \bar \lambda^\epsilon_2 \, .$$
                             \end{cor}
               
                \begin{proof}
 We observe that               
                $$\bar \lambda_2  ( B(0,1) \setminus \overline{B(0, \delta)})\geq \bar \lambda_2$$ 
                and that $\bar \lambda_2 ( B(0,1) \setminus \overline{ B(0, \delta)})$ is decreasing with respect to $\delta$. Hence,
                 $$\lim _{ \delta \rightarrow 0} \bar \lambda_2  ( B(0,1) \setminus \overline{B(0, \delta))}\geq \bar \lambda_2.$$ 
 In order to  prove the reverse inequality, we can use Theorem \ref{exp} jointly with Corollary \ref{convgammavar}, as well as Corollary \ref{delta}. Indeed,  for any  $\eta  >0$ let $\gamma_0<2$ such that,    for  all $ \gamma_0\leq \gamma<2$,  one has 
 $$\bar \lambda_2 \geq \bar \lambda_\gamma-\eta\, .$$ 
  Moreover,              let $\delta_0 =\delta( \gamma, \eta)$ be such that,  for any $\delta < \delta_0$, 
                $$\bar \lambda_\gamma  ( B(0,1) \setminus  \overline{B(0, \delta)})\leq \bar \lambda_\gamma + \eta\, .$$
It                 then follows
$$\bar \lambda_2 \geq \bar \lambda_\gamma  ( B(0,1) \setminus   \overline{B(0, \delta)})-2\eta \geq \bar \lambda_2( B(0,1) \setminus  \overline{B(0, \delta )})-2 \eta\, .$$
                The assertion concerning $\lim \bar \lambda^\epsilon_2 $ can be proved in the same way by using 
                Theorem \ref{th3}.
                 
                \end{proof}

In order to complete the proof of Theorem \ref{gamma=2}, it is enough to observe that statement (iii) immediately follows from statement (i), the definition of $\bar \lambda_2(F)$ and the ellipticity  inequalities \eqref{elliptic}.

                  \section{ The case $\gamma >2$ : Proof of Theorem  \ref{gamma>2} }
              
 This section is completely devoted to the proof of Theorem \ref{gamma>2}.     
 Let us assume, by contradiction, that for $\gamma>2$ there exists  $u\in C^2(B(0,1)\setminus \{0\})$  positive and radial, satisfying 
  $${\cal M}^-( D^2 u) \leq  -\mu u r^{-\gamma}\quad \hbox{ in } B(0,1)\setminus \{0\}\, , $$
  for some $\mu>0$. 
 
 Then, arguing  as in the proof of Lemma \ref{fabiana}, it follows that $u'(r)$ has constant sign in a right neighborhood of zero. If  $u^\prime (r)\geq 0$ for $r$ small, then, by the equation,  $u^{ \prime \prime}(r) \leq 0$ and then we  would have 
   $$(u^\prime r^{\tilde N_+-1})^\prime \leq -\frac{\mu u  r^{\tilde N_+-1-\gamma}}{\Lambda}\, .$$
 This implies that $u^\prime r^{\tilde N_+-1} $ has a nonnegative limit for $r\to 0$, and if this limit is strictly positive,  we get that $u$ becomes large negative as $r\to 0$, a contradiction.  Then the limit is zero, and then from the inequality above we get 
    $u^\prime(r) r^{\tilde N-1} <0$, a contradiction again. Therefore, we have 
    $u^\prime(r) \leq 0$ for $r>0$ sufficiently small. Then, we observe that, whatever is the sign of $u^{\prime \prime}$, one has 
    $$u^{\prime \prime} + (\tilde{N}_--1) \frac{u^\prime}{r} \leq -\frac{\mu}{\Lambda} u\,  r^{-\gamma}\, .$$
   Thus, 
        \begin{equation}\label{gamma2} ( u^\prime r^{ {\tilde N_--1}})^\prime \leq -\frac{\mu}{\Lambda} 
        u \, r^{\tilde{N}_--1-\gamma}< 0\, ,
     \end{equation}
 and then
      $u^\prime(r) r^{ {\tilde N_--1}}$ has a  non positive limit as $r\to 0$. If the limit is strictly negative, then   $u^\prime(r) r^{{\tilde N_--1}} \leq  -c <0$ in a neighborhood of zero.  
          Then,  we get
           $$u(r) \geq c_1 r^{2-{\tilde N_-}}$$
for $r$ small and a positive $c_1$. Hence, by integrating \eqref{gamma2} between $r$ and $s>r$ sufficiently small, we deduce
      $$ -u^\prime (s)s^{{\tilde N_--1}} + u^\prime (r) r^{{\tilde N_--1}} \geq c_2 \int_ r^s  t^{1-\gamma}  dt = c_2 
      \frac{s^{2-\gamma}- r^{2-\gamma}}{2-\gamma}$$
       and, since $\gamma >2$, this yields  $u^\prime(r)>0$ for $r$ small  enough:
        a contradiction.

 Thus, one has $\lim_{r\to 0} u^\prime r^{{\tilde N_--1}} = 0$ and, by \eqref{gamma2},   $u^\prime(r) <0$ for $r>0$. 
 
 Next, by an inductive argument analogous to the one used in the proof of Proposition \ref{signederiveeseconde}, we 
       prove that for all  integer $j\geq 0$ such  that $ {\tilde N_--2}+ j( 2-\gamma) >0$ and for $r$ sufficiently small, one has, for some $c_j>0$,
     \begin{equation}\label{eqj}u(r) \geq c_j r^{j(2-\gamma)}.
     \end{equation}
Indeed, \eqref{eqj} holds true for $j=0$, since $u^\prime<0$ and  $u$ is positive.   Let us suppose that (\ref{eqj}) is true for 
$j$   and   that ${\tilde N_--2}+ ( j+1)(2-\gamma) > 0$. Then, by (\ref{gamma2}),   
      $$ -u^\prime (r) r^{\tilde N_--1}  \geq \frac{\mu}{\Lambda} c_j \int_0^{r} s^{\tilde{N}_--1+ j ( 2-\gamma)-\gamma} ds =  \frac{\mu}{\Lambda} \frac{c_j}{\tilde{N}_--2+( j+1)( 2-\gamma)}  r^{{\tilde N_--2}+( j+1)  ( 2-\gamma)}\, ,$$
which, by integration, yields (\ref{eqj}) for $j+1$. 
            
 Now, let us assume  that $\frac{\tilde{N}_--2}{\gamma-2}$ is not  integer. Then, using estimate \eqref{eqj} with $j=\left[ \frac{\tilde{N}_--2}{\gamma-2}\right]$ jointly with \eqref{gamma2}, we deduce for $r_0>r>0$
 $$- u^\prime(r_0)  r_0^{{\tilde N_--1}}+ u^\prime (r) r^{{\tilde N_--1}}  \geq \frac{ \mu\, c_j}{\Lambda\, \left(\tilde{N}_--2+ (j+1)(2-\gamma)\right)}\left( s^{\tilde{N}_--2+ (j+1) (2-\gamma)}\right|_r^{r_0}\, .$$
       Since $\tilde{N}_--2+ (j+1) (2-\gamma)<0$, this yields the contradiction
$$\lim_{r\to 0} u^\prime (r) r^{{\tilde N_--1}}=+\infty\,.$$
On the other hand, if $ \frac{\tilde{N}_--2}{\gamma-2} = j+1$ is integer, then $\tilde{N}_--2+j(2-\gamma)=\gamma-2>0$, and from \eqref{eqj}  it follows that
            $$ u(r)\geq c_j r^{\gamma -\tilde{N}_-}\, ,$$ 
hence            
            $$ u(r) r^{ \tilde N_--1-\gamma} \geq c_j r^{-1}\, .$$
From \eqref{gamma2} we then deduce, for $0<r<r_0$,
            $$ -u^\prime (r_0) r_0^{ \tilde N_--1} + u^\prime (r) r^{ \tilde N_--1} \geq \frac{\mu\, c_j}{\Lambda} \left(  \ln r_0 - \ln r\right)$$
             and we reach, also in this case, the contradiction
$$\lim_{r\to 0} u^\prime (r) r^{{\tilde N_--1}}=+\infty\,.$$

\end{document}